\documentclass[11pt]{article}
\usepackage[a4paper, left=3cm, right=3cm, top=2cm]{geometry}
\usepackage[all]{xy}
\usepackage{graphicx}

\usepackage[latin1]{inputenc}
\usepackage[T1]{fontenc}

\reversemarginpar

\usepackage{amsmath}
\usepackage{amssymb}
\usepackage{latexsym}
\usepackage{amsfonts}
\usepackage{epsfig}
\usepackage{picinpar}
\usepackage[osf,sc]{mathpazo}
\usepackage{microtype}
\usepackage{textcomp}
\newcommand{\C}{{\bf C}}
\renewcommand{\P}{{\bf P}}
\newcommand{\Z}{{\bf Z}}
\newcommand{\Q}{{\bf Q}}
\renewcommand{\a}{\alpha}
\renewcommand{\b}{\beta}
\renewcommand{\c}{\gamma}
\begin{document}
\title{Calabi-Yau operators of degree two}
\author{Gert Almkvist\footnote{The first author passed away in December 2018.
Gert was a remarkable man and a great mathematician. It was a privilege to
have known and worked with him.}\; and Duco van Straten}
\maketitle

\begin{abstract} We show that the solutions to the equations defining
the so-called {\em Calabi-Yau condition} for fourth order operators of 
degree two defines a variety that consists of ten
irreducible components. These can be described completely in
parametric form, but only two of the components seem to 
admit arithmetically interesting operators. We include a description 
of the $69$ essentially distinct fourth order Calabi-Yau operators of 
degree two that are presently known to us. 
\end{abstract}

\section*{\Large \bf \em Introduction}
The hypergeometric series
\[\phi_0(x)=\sum_{n=0}^{\infty} \frac{(5n)!}{(n!)^5} x^n=1+120x+113400x^2+\ldots \in \Z[[x]]\]
became famous after in the paper \cite{COGP}  its dual interpretation was discovered: on the one hand the series encodes enumerative information on rational curves on the general quintic Calabi-Yau threefold in $\P^4$ and on the other hand
it can be identified as a normalised period of the mirror quintic Calabi-Yau threefold. 
The power series is the unique (up to scalar) series solution of the hypergeometric operator
\[L:=\theta^4-5^5x(\theta+\frac{1}{5})(\theta+\frac{2}{5})(\theta+\frac{3}{5})(\theta+\frac{4}{5}),\;\;\;\theta=x\frac{d}{dx} .\]
The operator $L$ is the first of the family of fourteen hypergeometric 
fourth order operators related to mirror symmetry for complete intersections in weighted projective spaces. These operators are all of the form
\[\theta^4-Nx(\theta+\alpha_1)(\theta+\alpha_2)(\theta+\alpha_3)(\theta+\alpha_4)\]
where the integer $N$ is chosen as to make the coefficients of the normalised
holomorphic solution integral. The properties of these operators have been 
well-studied from various points of view \cite{Alm1}, \cite{Villegas}, \cite{DM}, \cite{BT}, \cite{SV}, \cite{Singh}, \cite{Kon}.\\

The operator $L$ is also the first member of the ever growing list
of {\em Calabi-Yau operators} \cite{AESZ}, \cite{vS}, a notion that was 
introduced  by the first author and W. Zudilin in \cite{AZ} by abstracting 
the properties of $L$. Calabi-Yau operators in this sense are essentially 
self-adjoint fourth order fuchsian operators with a point of maximal 
unipotent monodromy for which strong integrality properties are supposed 
to hold. In a sense the term {\em Calabi-Yau operator} is somewhat of a
misnomer, as there exist families of Calabi-Yau varieties, whose Picard-Fuchs 
operator have no MUM-point, hence are not Calabi-Yau operators in the above sense of the word. For a recent account, see \cite{CynkvS}\\  

Let us consider a general $N$th order differential operator written 
{\em in $\theta$-form}:
\[L=P_0(\theta)+x P_1(\theta)+\ldots+\ldots+x^rP_r(\theta),\;\;\;\theta:=
x\frac{d}{dx},\]
where the $P_i$ are polynomials of degree $N$. We call the largest $r$ with 
$P_r \neq 0$ the {\em degree} of the operator. The differential equation
\[ L \phi=0,\;\;\; \phi(x)=\sum_{n=0}^{\infty} a_n x^n\]
translates into the recursion relation
\[ P_0(n)a_n+P_1(n-1)a_{n-1}+\ldots+P_r(n-r)a_{n-r}=0\]
on the coefficients $a_n$ of the series $\phi(x)$, so the degree of $L$ 
is equal to the {\em length} of this recursion.\\

The roots of $P_0(\theta)$ are the  {\em exponents} of the operator at $0$. 
By translation one can define exponents  of $P$ at each point of the Riemann 
sphere $\P^1$. For the point $\infty$ we have to apply inversion $x\mapsto 1/x$;
the exponents then are the roots of $P_r(-\theta)$. The {\em Riemann symbol} is a table that contains for each 
singular point the corresponding exponents, that encode information about the
local ramification of the solutions. Logarithmic terms may and usually do
appear at points with exponents with integral difference. This leads to Jordan
blocks in the local monodromy, but these are usually left implicit in the
Riemann symbol.\\
We will usually assume that $N=4$ and all exponents at $0$ vanish, $P_0(\theta)=\theta^4$. The most salient feature of this situation is that
at $0$ we have a canonical Frobenius basis of solution 
$\phi_0, \phi_1, \phi_2,\phi_3$, where $\phi_k$ contains terms with 
$\log^k(x)$. The local monodromy has a Jordan block of maximal size, so $0$
is a point of maximal unipotent monodromy: MUM.

The $14$ hypergeometric operators have degree $1$ and have a Riemann-symbol
of the form 
\[\left\{ \begin{array}{ccc}
0&1/N&\infty\\
\hline
0&0&\alpha_1\\
0&1&\alpha_2\\
0&1&\alpha_3\\
0&2&\alpha_4\\
\end{array}
\right\}\]
exhibiting a {\em conifold singularity} at the point $1/N$, where the exponents
are $0,1,1,2$ and the local monodromy around that point is a {\em symplectic reflection}; there is a Jordan block of size two. The exponents at $\infty$ are $\alpha_1 \le \alpha_2 \le \alpha_3 \le \alpha_4$ and add up in pairs to one: 
$$\alpha_1+\alpha_4=\alpha_2+\alpha_3=1.$$

Currently we know over $500$  Calabi-Yau operators with degrees running up to 
$40$,  but in no way do we expect the current list to be complete. An update 
of the AESZ-list is in preparation, \cite{ACvS}.\\

The first Calabi-Yau operator of degree $>1$ is \#15 in the AESZ-table and
appeared  in \cite{BvS}:
\[ \theta^4-3x(3\theta+1)(3\theta+2)(7\theta^2+7\theta+2)-72x^2(3\theta+1)(3\theta+2)(3\theta+4)(3\theta+5)\]
Its Riemann symbol is
\[\left\{ \begin{array}{cccc}
0&1/216&-1/27&\infty\\[1mm]
\hline
0&0&0&\frac{1}{3}\\[1mm]
0&1&1&\frac{2}{3}\\[1mm]
0&1&1&\frac{4}{3}\\[1mm]
0&2&2&\frac{5}{3}\\[1mm]
\end{array}
\right\}
\]

Another operator of this type \#25 of the AESZ-list. It arose from the  
quantum-cohomology of the Grassmanian $G(2,5)$. The complete interesection $X(1,2,2)$ of hypersurfaces of degrees $1,2,2$ in the Pl\"ucker embedding is a Calabi-Yau threefold with the following mirror-Picard-Fuchs operator: 
\[\theta^4  -4x (2\theta + 1)^2(11 \theta^2  + 11 \theta + 3) -16 x^2(2 \theta + 1)^2  (2 \theta + 3)^2 .\]
Its holomorphic solution is
\[ \phi_0(x)=\sum_{n=0}^{\infty} { 2n \choose n}^2A_n x^n,\]
where
\[A_n=\sum_{k=0}^n {n \choose k}^2{n+k \choose k}\]
are the {\em small Ap\'ery numbers} that were use in Ap\'ery's proof 
of the irrationality of $\zeta(2)$.\\
 
At present we know $69$ essentially distinct Calabi-Yau operators of order two, 
to which all further such degree two operators can be related by simple 
transformations. The operators for which the instanton numbers are $0$ are
obtained as $Sym^3$ of a second order operator and do not count as a proper
Calabi-Yau operator and do not appear in \cite{AESZ}\\

But contrary to the hypergeometric case, it appears that these degree two
operators do not all fall in a {\em single} family. First of all, an operator 
of degree two can have four singularities, with a Riemann symbol of the form

\[\left\{ \begin{array}{cccc}
0&a&b&\infty\\
\hline
0&0&0&\alpha_1\\
0&1&1&\alpha_2\\
0&1&1&\alpha_3\\
0&2&2&\alpha_4\\
\end{array}
\right\}
\]
exhibiting {\em two distinct conifold points}. But in many cases the points
$a$ and $b$ coalesce, producing an operator with only three singular points.

As one can observe in the above two examples, the exponents at infinity are
symmetrically centered around the value $1$ and this seemed to be
the case for all further examples we found. It was {\sc M. Bogner}
who first found an example of an Calabi-Yau operator for which this 
is {\em not} the case: it is the operator
\[{\theta}^{4}-x( 216\,{\theta}^{4}+396\,{\theta}^{3}+366\,{\theta}^{2}-168\,\theta+30)+  36\,x^2\left( 3\,\theta+2 \right) ^{2} \left( 6\,\theta+7 \right) ^{2}
\]

with Riemann symbol
\[\left\{ \begin{array}{cccc}
0&1/108&\infty\\
\hline
0&0&2/3\\
0&1/6&2/3\\
0&1&7/6\\
0&7/6&7/6\\
\end{array}
\right\}
\]
and later two more such operators were found. The existence of such operators 
was initially a surprise to us, as it appears that within the group of 
operators with three singular points further distinctions can be made.
It is this circumstance that led to this paper.\\

The structure of the paper is as follows. In the first section we recall 
some basic facts about the Calabi-Yau condition and formulate the main 
result of this paper. In the second section we describe how the result can 
be obtained from a series of simple computations. In the third section 
we summarise the properties of the operators corresponding to each of the 
components that we find. In a final section we review the list all degree 
two Calabi-Yau operators that are presently known to us. In two appendices 
we include some further monodromy and modular form information on these 
examples, together with some properties of the lower order operators used 
in many constructions.\\

\section{\Large \bf \em The Calabi-Yau condition} 
The {\em adjoint} $L^*$ of a differential operator  $L$
of the form
\[L =\frac{d^n}{dx^n}+a_{n-1} \frac{d^{n-1}}{dx^{n-1}}+\ldots+a_{1}\frac{d}{dx}+a_0 \in \Q(x)[\frac{d}{dx}]\]
is obtained by reading the operator {\em backwards with alternating signs}:
\[L^* =\frac{d^n}{dx^n}-\frac{d^{n-1}}{dx^{n-1}}a_{n-1}+\ldots+\pm\frac{d}{dx}a_1 +(-1)^n a_0 \in \Q(x)[\frac{d}{dx}],\]
so that
\[L^*(y) =y^{(n)}- (a_{n-1}y)^{(n-1)}+\ldots+(-1)^n a_0 y\]

The operator $\mathcal{P}$ is called {\em essentially self-adjoint}, if 
there  exists a function  $\alpha \neq 0$ (in an extension of $\Q(x)$), 
such that
\[ L \alpha =\alpha L ^{*} .\]
It is easy to see that any such  $\alpha$ has to satisfy
the first order differential equation
\[ \alpha' = -\frac{2}{n} a_{n-1} \alpha .\]
This essentially adjointness is equivalent to the existence of 
an invariant pairing on the solution space, symmetric
if $n$ is odd and alternating if $n$ is even, \cite{Bogner1}, \cite{Bogner2}.
The condition of essential self-adjointness can be expressed by the vanishing 
of certain differential polynomials in the $a_i$ described in
\cite{AZ}, \cite{Bogner1}, \cite{Bogner2}, \cite{ABG}, called the 
{\em Calabi-Yau condition}.\\

For an operator of order two the essentially self-adjontness does not impose
{any conditions}, whereas for an differential equation of order three of the 
form
\[ y'''+a_2 y''+a_1y'+a_0y \]
the Calabi-Yau condition comes down to the vanishing of the quantity
\[\mathcal{W}:=\frac{1}{3}a_2''+\frac{2}{3}a_2a_2'+\frac{4}{27}a_2^3+2a_0-\frac{2}{3}a_1a_2-a_1'\;,\]
whereas for a fourth order operator
\[ y^{(iv)}+a_3 y'''+a_2 y''+a_1 y' +a_0 y = 0 \]
the quantity
\[\mathcal{Q}:= \frac{1}{2} a_2 a_3-a_1-\frac{1}{8}a_3^3+a_2'-\frac{3}{4}a_3 (a_3)'-\frac{1}{2} a_3''\]
has to vanish. In that case the associated differential equation satisified 
by the $2\times 2$  Wronskians reduces its order from $6$ to $5$.\\

For operators of higher order, one finds that more than one differential 
polynomial condition has to be satisfied; here we will be mainly concerned 
with operators of order four.\\

It is rather easy to fulfill the condition $\mathcal{Q}=0$.
For example, there exists one particularly nice family of operators 
of degree $n$, that we call the {\em main family}.\\

\centerline{\Large \bf \em Proposition} 
\vskip 10pt
{\em  If $P,Q,R$ are polynomials of degree $n$, then the operator
\[L= \theta^2 P \theta^2+\theta Q \theta +R\]
\[  = P\theta^4+2 P'\theta^3+(Q+P'')\theta^2+Q'\theta+R\]
has degree $n$ and satisfies the Calabi-Yau condition.
If the constant term of $P$ is $1$ and of $Q$ and $R$ is zero, then the
exponents of $L$ at $0$ are all zero. 
If the roots of $P$ are all distinct, $L$ has $n$ singular points
with exponents $0,1,1,2$ and a further singular point at infinity.}
(Here $P':=\theta(P)=x\partial P/\partial x$, etc.)

\vskip 10pt

In this paper we take a closer look at differential operators of the form 
\[\theta^4+x(a\theta^4+b\theta^3+c\theta^2+d\theta+e)+fx^2(\theta+\alpha)(\theta+\beta)(\theta+\gamma)(\theta+\delta)\]
If $f \neq 0$, then the numbers $\alpha,\beta,\gamma,\delta$ are then the 
characteristic exponents of the operator at $\infty$.\\

For this operator the quantity $\mathcal{Q}$ is a rational function with 
\[ x^2(1+ax+fx^2)^3\]
as denominator; the singular points of the operator are the roots of
this polynomial, together with the point $\infty$. \\

The numerator is a polynomial $Q$ of degree $5$ in the variable $x$, so we
can write
\[ Q=Q_0+Q_1x+Q_2x^2+Q_3x^3+Q_4x^4+Q_5x^5, \]
where the coefficients $Q_i$ are rather complicated polynomials in the 
parameters $a,b,c,d,e,f$ and the exponents $\alpha,\beta,\gamma, \delta$:
\[ Q_0,\;\;Q_1,\;\;Q_2,\;\;Q_3,\;\;Q_4,\;\;Q_5 \in \Q[a,b,c,d,e,f,\a,\b,\c,\delta].\]
So the Calabi-Yau condition is equivalent to  six polynomial equations 
\[Q_0=Q_1=Q_2=Q_3=Q_4=Q_5=0\]
which define a certain {\em affine algebraic set}
\[ X=V(Q_0,Q_1,\ldots,Q_5) \subset {\bf C}^{10} .\]

We will give a complete description of the set $X$, in particular
we determine the irreducible components of $X$. As we have six equations in 
$10$ variables, a first trivial remark is that all irreducible components of 
$X$ have dimension {\em at least four}. We will give explicit parametrisations 
of all the irreducible components. Note also that the polynomials $Q_i$ do not
depend on the accessory parameter $e$, so the solution set $X$ contains the
$e$-line as trivial factor.\\

\centerline{\Large \bf \em Theorem}
\vskip 10pt 
{\em The algebraic set $X$ defined by the condition $Q(a,b,c,d,e,f,\alpha,\beta,\gamma,\delta)=0$
is the union of ten irreducible components, of which there are seven 
with $f\neq 0$.}\\

Note that if $f=0$ we are dealing with an operator of lower degree, so 
will from now on assume that $f \neq 0$ and we will ignore the three 
components with $f=0$. The seven remaining components differ in the 
behaviour of the average exponents at infinity:

\[ \sigma:=\frac{\alpha+\beta+\gamma+\delta}{4} .\]

There is a single component for which $\sigma$ is not fixed, that we will 
call the {\em transverse component}. On each of the remaining six components
$\sigma$ has a fixed value. The dimensions and values of $\sigma$ that 
appear are summarised in the following diagram.
\[
\begin{array}{|c||c|c|c|c|c|}
\hline
\sigma&0&1&2&3&4\\
\hline
      & &   &6 &   &\\ 
\dim  & & 5 &  &  5& \\
      &4&   &4  &  &4\\
\hline
\end{array}
\]
The two five dimensional components are related by a simple transformation, 
as are the three four dimensional components.\\

All Calabi-Yau operators known to us, belong to the six-dimensional 
{\em main component $M$:}

\[ M: \theta^4+x(a\theta^4+2 a \theta^3+(a+d)\theta^2+d\theta+e)+fx^2(\theta+\a)(\theta+\b)(\theta+2-\a)(\theta+2-\b) \]

or the  {\em transverse component $T$} that we describe later.

We remark furthermore that the ideal 
\[I=(Q_0,Q_1,Q_2,Q_3,Q_4,Q_5) \subset \C[a,b,c,d,e,f,\alpha,\beta,\gamma,
\delta]\]
is, from an algebraic point of view, rather complicated. It is not radical 
and the components appear with non-trivial multiplicities in the primary 
decomposition.\\

For third order operators of degree two
\[ \theta^3+x(a\theta^3+b\theta^2+c\theta+d)+f x^2(\theta+\a)(\theta+\b)(\theta+\c)\]
the Calabi-Yau condition $\mathcal{W}=0$ leads in a similar way to an
algebraic set
\[ Y \subset \C^8\]
defined by a system of polynomial equations
\[ W_0=W_1=W_2=W_3=W_4=W_5=0\]
where now the polynomials
\[ W_0,\;\;W_1,\;\;W_2,\;\;W_3,\;\;W_4,\;\;W_5 \in \Q[a,b,c,d,f,\a,\b,\c].\]
are in eight variables. The analysis of this ideal is very similar to the
ideal for the fourth order operators, but hardly simpler. The result is also 
very similar to the case of order four:\\

\centerline{\Large \bf \em Theorem}
\vskip 10pt 
{\em The algebraic set $Y \subset \C^8$ defined by the condition 
$W(a,b,c,d,f,\alpha,\beta,\gamma)=0$
is the union of ten irreducible components, of which there are seven with 
$f\neq 0$.
}\\

The seven remaining components differ in the behaviour of the (doubled) average
exponent
\[ \sigma=2 \frac{\a+\b+\c}{3} .\]
Again there is a single {\em transverse component} on  which the value of 
$\sigma$ is not fixed. On each of the remaining six components
$\sigma$ has a fixed value. The dimensions and values of $\sigma$ that 
appear are summarised in the following diagram.

\[
\begin{array}{|c||c|c|c|c|c|}
\hline
\sigma&0&1&2&3&4\\ 
\hline
       &  &  & 4 &   &\\ 
\dim   &  & 3&   & 3 &\\
       &2 &  & 2 &   & 2\\
\hline
\end{array}
\]

The main component $(m)$ is the one with $\sigma=2$, $\dim=4$ and is given by
\[(m):  \theta^3+x(2\theta+1)((c-2d) \theta^2+(c-2d)\theta+ d)+f x^2(\theta+\a) (\theta+1) (\theta+2-\alpha),  \]
whereas on the transverse component $(t)$ we have the operator
\[\theta^3+ \frac{a}{4} x(4\theta^3+3(\alpha+\beta)\theta^2+(\a+\b+2\a \b)\theta +\alpha \beta)+ (\frac{a}{4})^2 x^2 (\theta+\a)(\theta+\frac{\a+\b}{2})(\theta+\b)
\]
Again, the operators on the two three dimensional components are related by a
simple transformation, as are the three two dimensional components.

One may speculate that the above patterns extends to higher order operators, 
but we will not try to pursue it here.\\

\section{\Large \bf \em Proof of the theorem}
The proof of the above theorem is purely computational and consists of 
systematically solving of the equations using a computer algebra system 
like {\sc Maple}. The miracle here is that this procedure actually works. 
We will  describe the most important steps in the process.\\

\subsection{\bf \em Preliminary reductions}

{\em The value of $b$:} The polynomial $Q_0$ is found to be
\[-4 b - 8 d + 8 c\]
and this allows us directly to eliminate the variable $b$ and thus we will put
from now on
\[ b=2(c-d)\;\;.\]

{\em Symmetry of the exponents:} The polynomial $Q_5$ factors nicely as
\[-f^3(\alpha+\beta-\delta-\epsilon)(\alpha+\delta-\beta-\epsilon)(\alpha+\epsilon-\beta-\delta)\]

As we assume that $ f \neq 0$, this only vanishes if one of the three other
factors vanish, meaning that {\em the exponents come in two pairs that add up to the same value} that we shall call $\sigma$. So without loss of generality we can assume to have the pairs $\alpha, \sigma-\alpha$ and $\beta, \sigma-\beta$ and we will take the operator to be of the form
\[ \theta^4+x(a\theta^4+b\theta^3+c\theta^2+d\theta+e)+fx^2(\theta+\alpha)(\theta+\beta)(\theta+\sigma-\alpha) (\theta+\sigma-\beta), \]
with $b=2(c-d)$ as our starting point.\\

After solving the equations $\Q_0=0,\Q_5=0$, one is left to solve the polynomial equations $Q_1=Q_2=Q_3=Q_4=0$. 
For sake of concreteness, we give the explicit form of the polynomials:

\[{Q_1}:=8\,f{{ \alpha}}^{2}{ \sigma}+16\,f{ \sigma}\,{ \beta}-8\,f{{\sigma}}^{2}{ \beta}-16\,{c}^{2}+40\,cd-24\,{d}^{2}-8\,f{ \alpha}\,{{ \sigma}}^{2}-16\,f{{ \alpha}}^{2}+\]
\[+16\,f{{ \sigma}}^{2}-16\,f{{ \beta}}^{2}-32\,f{ \sigma}-24\,da+16\,ca+16\,f{ \alpha}\,{ \sigma}+8\,f{ \sigma}\,{
{ \beta}}^{2},\]

\[ {Q_2}:=-8\,f{{ \alpha}}^{2}c-8\,f{{ \beta}}^{2}c+8\,f{{ \beta}}^{2}d+
8\,f{{ \alpha}}^{2}d-8\,cad+8\,f{{ \sigma}}^{2}c+72\,f{ \sigma}\,d-
8\,f{{ \sigma}}^{2}d+\]
\[+24\,{c}^{2}d-24\,c{d}^{2}+8\,{c}^{2}a-64\,cf{ \sigma}+24\,f{ \alpha}\,{ \sigma}\,a
+24\,f{ \sigma}\,{ \beta}\,a-16\,f{{ \sigma}}^{2}{ \beta}\,a+16\,f{{ \alpha}}^{2}{ \sigma}\,a+\]
\[-16\,f{ \alpha}\,{{ \sigma}}^{2}a+16\,f{ \sigma}\,{{ \beta}}^{2}a
-8\,f{ \alpha}\,{ \sigma}\,d+8\,f{\alpha}\,{ \sigma}\,c-8\,{c}^{3}+8\,{d}^{3}+
8\,f{ \sigma}\,{ \beta}\,c-8\,f{ \sigma}\,{ \beta}\,d+\]
\[-64\,df-8\,d{a}^{2}-24\,f{ \sigma}\,a-24\,f{{\alpha}}^{2}a+48\,cf+24\,f{{ \sigma}}^{2}a-24\,f{{ \beta}}^{2}a,\]

\[{Q_3}=8\,f \left( 2\,f{ \sigma}\,{ \beta}+2\,f{ \alpha}\,{ \sigma}+f{ \sigma}\,{{ \beta}}^{2}-f{{ \sigma}}^{2}{ \beta}+f{{ \alpha}}^{2}{\sigma}-f{ \alpha}\,{{ \sigma}}^{2}-2\,f{{ \beta}}^{2}+\right.\]
\[-3\,da+{{ \sigma}}^{2}{a}^{2}+f{{ \sigma}}^{3}-{{ \beta}}^{2}{a}^{2}-{ \sigma}\,{a}^{2}-{ \sigma}\,{ \beta}\,ad+{ \sigma}\,{ \beta}\,ac-{ \alpha}\,{ \sigma}\,ad+{ \alpha}\,{ \sigma}\,ac+\]
\[-{{ \alpha}}^{2}{a}^{2}-2\,c{ \sigma}\,a-2\,f{{ \alpha}}^{2}-4\,f{{ \sigma}}^{2}+4\,f{ \sigma}+4\,{c}^{2}-7\,cd+3\,{d}^{2}
-{ \alpha}\,{{ \sigma}}^{2}{a}^{2}+\]
\[+{ \sigma}\,{{ \beta}}^{2}{a}^{2}+6\,{ \sigma}\,cd+3\,{ \sigma}\,ad-{{ \beta}}^{2}ac
+{{ \beta}}^{2}ad-{{ \alpha}}^{2}ac+{{ \alpha}}^{2}ad+{{ \sigma}}^{2}ac+\]
\[\left.-{{ \sigma}}^{2}ad-3\,{ \sigma}\,{d}^{2}-3\,{ \sigma}\,{c}^{2}+{ \alpha}\,{ \sigma}\,{a}^{2}+
{ \sigma}\,{ \beta}\,{a}^{2}-{{ \sigma}}^{2}{ \beta}\,{a}^{2}+{{ \alpha}}^{2}{ \sigma}\,{a}^{2} \right), 
\]

\[
{ Q_4}=8\,{f}^{2} \left( -2\,c-{{ \alpha}}^{2}c-2\,{{ \sigma}}^{2}c-3\,{ \sigma}\,d+2\,{{ \sigma}}^{2}d+{{ \beta}}^{2}d-2\,{{ \sigma}}^{2}a-{
{ \beta}}^{2}a-{ \alpha}\,{ \sigma}\,d+\right.\]
\[{ \alpha}\,{ \sigma}\,c+{ \sigma}\,{
 \beta}\,c-{ \sigma}\,{ \beta}\,d+{ \sigma}\,{ \beta}\,a-{ \alpha}\,{{
 \sigma}}^{2}a-{{ \sigma}}^{2}{ \beta}\,a+{ \alpha}\,{ \sigma}\,a+{{ \alpha
}}^{2}{ \sigma}\,a+{{ \beta}}^{2}{ \sigma}\,a+\]
\[\left.{ \sigma}\,a-{{ \alpha}}^{2}
a+{{ \sigma}}^{3}a+4\,c{ \sigma}-{{ \beta}}^{2}c+{{\alpha}}^{2}d
 \right) .\]

Somewhat to our surprise, it turns out to be rather easy to give a complete
solution to these equations. It will be useful to introduce the parameters
\[ A:=\alpha(\sigma-\alpha),\;\;B:=\beta(\sigma-\beta)\]
and use these instead of $\alpha$ and $\beta$. 
One has
\[(\theta+\alpha)(\theta+\sigma-\alpha)(\theta+\beta)(\theta+\sigma-\beta)=(\theta^2+\sigma \theta +A)(\theta^2+\sigma \theta+B)\]
\[=\theta^4+2\sigma \theta^3+(\sigma^2+(A+B))\theta^2+\sigma(A+B)\theta+AB .\]
We will see that most expressions only depend on 
\[\Delta:=-(A+B)=\alpha^2+\beta^2-\sigma(\alpha+\beta) .\]

\subsection{\em The polynomial $Q_4$}

The polynomial $Q_4$ is an expression that is $f^2$ times an expression linear in the variable $a$. 
As we will only consider cases with $f \neq 0$, we  can solve for $a$ and 
find
\[
a=\frac{(2(\sigma-1)^2+\Delta)c-(\sigma(2\sigma-3)+\Delta)d}{(\sigma-1)(\sigma^2-\sigma+\Delta)}\;\;,
\]
Note that this value for $a$ introduces the denominator
\[ (\sigma-1)(\sigma(\sigma-1)+\Delta) .\]
So if we work further with this operator, we are implicitly assuming that
these two factors do not vanish. If one of them does vanish, we should go back to 
the previous step, impose these conditions and compute further. For now, we will assume that these factors are non-zero and come back to these cases later.\\

\subsection{\bf \em The case $\sigma=2$}
It appears that with this value of $a$ the quantity $\mathcal{Q}$ becomes
divisible by $\sigma-2$. So in this remarkable $\sigma=2$ case
we obtain an operator family with $\mathcal{Q}=0$.
The value of $a$ simplifies to
\[a=\frac{(2+\Delta)c-(2+\Delta)c}{1\cdot (2+\Delta)}=c-d\;,\]
and one obtains what we call the {\em main component:}
\[P:=\theta^4+x((c-d)\theta^4+2(c-d)\theta^3+c\theta^2+d\theta+e)+\]
\[ +fx^2(\theta+\alpha)(\theta+\beta)(\theta+2-\alpha) (\theta+2-\beta)\;,\]
depending $c,d,e,f$ and the exponents $\alpha$ and $\beta$ as free parameters.
Almost all operators of degree  two from the AESZ-list are of this type.
We obtain our first component of dimension $6$.


\subsection{\bf \em The case $\sigma \neq 2$}

If however $\sigma \neq 2$ we have to analyse further. Still using the above value for $a$, the coefficient of $x$ of $\mathcal{Q}$ is a 
complicated expression in $c,d,e,f$ and the exponents, that however is factored by  {\sc Maple} instantly into four factors.
First, there are the factors 
\[ (\sigma-1),\; (\sigma-2),\;(\sigma (\sigma-1)+\Delta) .\]\\ 
We assume first that these are non-zero, so one is forced put to zero the 
fourth factor, which leads to a specific value for $f$:
\[
f=(2c-3d){\frac { ((\sigma-1)^2+\Delta)c-(\sigma(\sigma-2)+\Delta)c}{ \left( \sigma-1 \right)  \left (2 \sigma +\Delta\right)  \left( {{\it \sigma}}^{2}-{\it \sigma} +\Delta \right)}} .
\]

Here a new factor $2\sigma+\Delta$ is introduced in the denominator, so we
have to come back later to the case $2\sigma+\Delta=0$.


\subsubsection{\bf \em The subcase $\sigma=3$}
Using the above values for $a$ and $f$, the quantity $\mathcal{Q}$ only contains terms with $x^2$ and $x^3$. It factors out the factor $\sigma-3$, which leads to another remarkable operator family
with $\mathcal{Q}=0$:
\[a =\frac{(8+\Delta)c-(9+\Delta)d}{2(6+\Delta)},\;\;\;f = (2c-3d) \frac{(4+\Delta)c-(3+\Delta)d}{2(6+\Delta)^2} .\]
It has $c,d,e$ and $\alpha,\beta$ as free parameters. This makes up a second component, of dimension $5$.                                          

\subsubsection{\bf \em The subcase $\sigma \neq 3$}
Still with the given values for $a$ and $f$, if $\sigma \neq 3$, one has to analyse the quantity $\mathcal{Q}$ further.  {\sc Maple} manages to factor the coefficient of $x^2$  into a product of eight factors, of only three were not supposed to be non-zero at this stage. Each of these factors leads to a linear dependence between $c$ and $d$.

{\bf Case A:}
\[c=(\Delta-\sigma)\frac{d}{\Delta} .\]
Upon substitution we obtain an operator that satisfies the condition $\mathcal{Q}=0$. For this operator one has
\[a=-2 \frac{d}{\Delta},\;\;\;b=-2\sigma \frac{d}{\Delta},\;\;\;c=-(\sigma-\Delta) \frac{d}{\Delta},\;\;\;f= (\frac{d}{\Delta})^2,\]
and has $d, e$ and the exponents $\alpha, \beta,\sigma$ as free parameters. So we obtain a further irreducible component of dimension $5$.
The remarkable thing is that now we have only $\Delta$ apearing in the 
denominator.\\

The operator of {\sc Bogner} mentioned in the introduction is an instance of 
case A for the  parameter choice:
\[\alpha=\beta=2/3,\;\;\sigma=11/6,\;\;d=-168,\;\; e=30 .\]

{\bf Case B:}

\[c=d \frac{\sigma^2-2\sigma+\Delta}{(\sigma-1)^2+\Delta} .\]


In this case $f$ reduces to $0$, so we do not get an operator of degree two.\\

Case C:

\[c=d \frac{\sigma^2-6\sigma+4+\Delta}{(\sigma-1)(\sigma-4)+\Delta}\]

In this case the condition $\mathcal{Q}=0$  reduces to special relations between 
the exponents. The exceptional cases are:
\[\sigma \in \{0,1,2,3,4\}\]
or
\[\left( {{\it \sigma}}^{2}-2\sigma+4+\Delta \right)=0,\;\; \left( {{\it \sigma}}^{2}-3\sigma+4+\Delta \right)  =0,\;\;  \left( {{\it \sigma}}^{2}-5\sigma +4+\Delta \right)  =0\]

This concludes the list of all possibilities. We still have to backtrack 
some of the cases:
 
\subsection{\bf \em Backtracking the remaining cases}

Earlier we use the value for $a$ which involved the denominator
\[(\sigma-1)(\sigma^2-\sigma+\Delta)\;.\] 
The introduction of $f$ led to a further denominator 
\[2\sigma+\Delta\] 
The analysis of Case A gave $\Delta$ as denominator, and Case C gave 
further factor  $(\sigma-1)(\sigma-4)+\Delta$ in the denominator.
In each of the cases we have to backtrack and see if we find additional
solutions.\\

\subsubsection{\bf \em The case $\sigma=1$}
Looking at the operator with $\sigma=1$, we find from the
coefficient of $x$ of $\mathcal{Q}$ the value of $f$. Using this value
factoring the coefficient of $x^2$ of $\mathcal{Q}$ we find there is a unique
value of $c$:
\[c=d\frac{\Delta-1}{\Delta} .\]
Here a denominator  $\Delta$ appears, so we assume this to be $\neq 0$;
but the case $\Delta=0$ has to be analysed anyway.

\subsubsection{\bf \em The case $\Delta=0$}
To analyse the further cases, it is useful to write the operator in terms
of $A$ and $B$ instead of $\alpha$ and $\beta$, so that $\Delta=C$ can
be implemented by putting $B=C-A$. We start with $C=0$.
The following seven components make up the intersection with $\Delta=0$:\\

Case 1) $\;\;\sigma=\Delta=0,\;\;b,\;c,\;d=0$:
\[\theta^4+x(a\theta^4+e)+fx^2(\theta^4+u).\]

Case 2) $\;\;\sigma=1,\;\;\Delta=0,\;\;d=0,\;\;f=(a-c)c$;
\[\theta^4+x(a\theta^4+2c\theta^3+c\theta^2+e)+(a-c)c(\theta^4+2\theta^2+3+\theta^2+u)\]

Case 3) $\;\;\sigma=2,\;\;\Delta=0,\;\;a=c-d$. This is inside the  main component.\\

Case 4) $\;\;\sigma=2,\;\; \Delta=0,\;\;a=d,\;\;c=2d.\;\;f=d^2/4$.\\

Case 5) $\;\;\sigma=3,\;\; \Delta=0,\;\;a=8c-9d/12,\;\;f=(2c-3d)(4c-3d)/72$\\
This is inside the  $\sigma=3$ component.\\

Case 6) $\;\;\sigma=4,\Delta=0,\;\;d=0,\;\;c=2a,\;\;f=a^2/4$.\\

Case 7) $\;\;d=\Delta=0$, there is the family with parameters $c,e,\sigma,\alpha$
\[a=\frac{2c}{\sigma},\;\;b=2c,\;\;d=0,\;\;f=\frac{a^2}{4}.\]

\subsubsection{\bf \em Remaining cases}
There are more cases to check and this could be done along the
above sketched lines. Instead of this, we used the computer
algebra system {\tt Singular}. There is a library called
{\tt primdec.lib} containing the algorithms for doing primary
decomposition. On our computer, {\tt Singular} could not straight away compute the
primary decomposition of the ideal $\mathcal{I}$ defined by the 
vanishing of the $x$-coefficients of $\mathcal{Q}$. But after
taking the ideal quotient of $\mathcal{I}$ by the ideals already
found, it turned out to be able to find the last component.
In total we thus found seven irreducible components of the variety
defined by $\mathcal{Q}=0$. 

\subsection{\bf \em Summary of the seven components}

In total we have found seven components. Membership to a component
is defined by certain polynomial relations between the coefficients.
Below we give the defining relations for each of the components, the
operator with its coefficients together with its generic Riemann symbol.\\ 

{\bf \em The big $\sigma=2$ component}, also called the {\bf \em main component}
is defined by the linear equations

\[\sigma-2=0,\;\;a-(c-d)=0.\]

For a general member of this component the Riemann symbol is:
\[
\left\{
\begin{array}{cccc}
0&*&*&\infty\\ 
\hline
0&0&0&\alpha\\
0&1&1&\beta\\
0&1&1&2-\beta\\
0&2&2&2-\alpha
\end{array}
\right\}
\] 
where $*$ are the two roots of $1+ax+fx^2=0$.
The operator depends on the six parameters $a,d,e,f,\alpha,\beta$.
The parameters $a,f$ determine the position of the singular fibres.
The parameters $\alpha,\beta$ determine the four exponents at $\infty$ which
are symmetric around the value $2$. The parameters $d,e$ are accessory 
parameters.\\
However, if the two roots coincide, $a^2=4f$, then we obtain an operator
with three singular points.\\ 

{\bf \em The $\sigma=1$ component} is defined by the four conditions:
\[\sigma-1=0,\;\;c\Delta+d(\Delta-1)=0,\;\;\]
\[d(a-(c-d))+f\Delta=0,\;\;c(a-(c-d))+f(\Delta-1)=0\]

The Riemann symbol of the operator is:
\[
\left\{
\begin{array}{cccc}
0&x&y&\infty\\
\hline
0&0&0&\alpha\\
0&1&1&\beta\\
0&1&2&1-\beta\\
0&2&3&1-\alpha\\
\end{array}
\right\}
\]
where the special points $x$ and $y$ are located at:
\[ x=\frac{\Delta}{d},\;\;\; y=\frac{\Delta}{d+a\Delta}.\]
The point $y$ is not really a singular point of the differential
operator, as the exponents $0,\;1,\;2,\;3$ are those of a regular point
of the differential equation. Indeed, as it can be checked that no 
logarithms occur. \\

{\bf \em The $\sigma=3$ component} is defined by the conditions:
\[
\begin{array}{rcl}
0&=&\sigma-3\\
0&=&4a^2-4c^2-3ad+11cd-7d^2-3fA+20f\\
0&=&2aA-cA+dA-12a+8c-9d\\
0&=&2ac-2c^2-3a^2aA-cA+dA-12a+8c-9d\\
0&=&2c^2A-5cdA+3d^2A+2fA^2-8c^2+18cd-9d^2-24fA+72f\\
\end{array}
\]

The Riemann symbol of the operator is:
\[
\left\{
\begin{array}{cccc}
0&x&y&\infty\\
\hline
0&0&-1&\alpha\\
0&1&0&\beta\\
0&1&1&3-\beta\\
0&2&2&3-\alpha\\
\end{array}
\right\}
\]

where the singular points are located at:
\[x=\frac{6+\Delta}{2c-3d},\;\;y=\frac{-2(6+\Delta)}{4c-3d+(c-d)\Delta} .\]
The operator has five parameters $c,d,e,\alpha,\beta$.
$\alpha, \beta$ are parameters that determine the exponents at infinity,
symmetric around $3$. The parameters $c,d$ determine, together with 
$\alpha$ and $\beta$ the position of the two singularities. 
$e$ is a single accessory parameter.\\

{\bf \em  The $\sigma=0$ component} is defined by the conditions:
\[\sigma=0,\;\;c=0,\;\;d=0,\;\;\Delta=0.\]

The generic Riemann symbol for this operator is:
\[
\left\{
\begin{array}{cccc}
0&x&y&\infty\\
\hline
0&0&0&\alpha\\
0&1&1&i\alpha\\
0&2&2&-\alpha\\
0&3&3&-i\alpha\\
\end{array}
\right\}
\]

where the special points $x$, $y$ are located at the solutions of 
$1+ax+fx^2=0$.
The operator depends on four parameters $a,e,f,\alpha$.\\

{\bf \em The small  $\sigma=2$ component} is defined by the four conditions:

\[\sigma-2=0,\;\;4ad-d^2-16f=0,\;\;2c-3d=0,\;\;\Delta+2=0 .\]

The Riemann symbol for a general member of this component is:
\[
\left\{
\begin{array}{cccc}
0&x&y&\infty\\
\hline
0&-1&0&2+\alpha\\
0& 0&1&2+i\alpha\\
0& 1&2&2-\alpha\\
0& 2&3&2-i\alpha\\
\end{array}
\right\}
\]
where the special points are 
\[x=-\frac{4}{d},\;\;y=-\frac{4}{4a-d} .\]

{\bf  \em The $\sigma=4$ component} is defined by the four conditions:
\[\sigma-4=0,\;\;6a-c=0,\;\;2c-3d,\;\;\Delta+8=0.\]
The Riemann symbol of this operator is:
\[
\left\{
\begin{array}{cccc}
0&x&y&\infty\\
\hline
0&-1&-1&4+\alpha\\
0& 0&0&4+i\alpha\\
0& 1&1&4-\alpha\\
0& 2&2&4-i\alpha\\
\end{array}
\right\}
\]
where the singular points $x,y$ are solutions to
$1+ax+fx^2=0$.\\

{\bf The transverse component} finally is defined by the equations
\[
\begin{array}{rcl}
0&=&a^2-4f,\\
0&=&a\sigma-2(c-d),\\
0&=&a\Delta+2d,\\
0&=&-f\Delta^2+d^2,\\
0&=&ad+2f\Delta,\\
0&=&d\sigma-c\Delta+d\Delta,\\
0&=&f\sigma\Delta-f\Delta^2+cd,\\
0&=&-f\sigma^2+f\sigma\Delta+c^2-cd,\\
0&=&ac-2f\sigma+2f\Delta.\\
\end{array}
\]

Solving for $a,b,c,f$ we find the family of operators
\[T(d,e,\alpha,\beta,\sigma):= \theta^4+\]
\[-x \frac{1}{\Delta}\left(2d\theta^4+2d \sigma \theta^3+(\sigma-\Delta)\theta)\right)+x(\theta+e)+\] 
\[-\left(\frac{d x}{\Delta}\right)^2(\theta+\alpha)(\theta+\beta)(\theta+\sigma-\beta)(\theta+\sigma-\beta)\]

The Riemann symbol of the generic member in this family is:
\[
\left\{
\begin{array}{cccc}
0&x&\infty\\
\hline
0& 0&\alpha\\
0& 1&\beta\\
0& 2-\sigma&\sigma-\beta\\
0& 3-\sigma&\sigma-\alpha\\
\end{array}
\right\}
\]
where $x=\frac{\Delta}{d}$.\\

The big and small $\sigma=2$ components do intersect:
\[\sigma=2,\;\;\Delta=-2,\;\;2a=d,\;2c=3d,\;f=d^2/16.\]
The transverse component intersects all other components.
No other pair of components can intersect, as they have 
a different $\sigma$-value.\\

\subsection{\bf The small components}
It turns out that the components with $\sigma=0,\;\;2,\;\;4$ are
closely related to each other and the solutions of the corresponding
differential equations can be related in a simple manner.
The families of operators in question are:\\

\[ 
\begin{array}{rcl}
P_0(a,e,f,A)&:=&\theta^4+x(a\theta^4+e)+x^2f(\theta^4-A^2)\\[3mm]
&&\\
P_2(a,d,e,A)&:=&\theta^4+x(a\theta^4+d\theta^3+(3/2)d\theta^2+d\theta+e)+\\[3mm]
&&+\frac{d(4a-d)}{16}x^2((\theta+1)^4-(A-1)^2)\\[3mm]
&&\\
P_4(a,e,f,A)&:=&\theta^4+x(a\theta^4+4a\theta^3+6a\theta^2+4a\theta+e)+\\[3mm]
&&+x^2 f(\theta^4+8\theta^3+24\theta^2+32\theta+A(8-A))\\[3mm]
\end{array}
\]

Their relation is most easily understood from the perspective of the operator
$P_2(a,d,e,A)$ and its exponents. The discriminant $\Delta(x)$ of the operator $P_2(a,d,e,A)$ factors as $(1+dx/4)(1+(4a-d)x/4)$.
The exponents at the first factor are $-1, 0, 1, 2$ and at the second factor
$0, 1, 2, 3$. The exponents at infinity are of the form
\[ 1+\alpha,\;\; 1-\alpha,\;\; 1+i \alpha,\;\;1-i \alpha . \]
Multiplication of a solution to $P_2(a,d,e,A)$ by $(1+dx/4)$ will shift the exponents  at the first singular point by $+1$ and the exponents at infinity by $-1$, which then leads to an operator with two points with exponents 
$0,1,2,3$ and 
\[ \alpha, -\alpha, i \alpha,-i \alpha \]
as exponents at infinity, i.e with $\sigma=0$ and $\Delta =0$.
Similarly, division of a solution to $P_2(a,d,e,A)$ by $(1+(4a-d)x/4)$ leads to a shift at the second singularity by $-1$ and hence produces two 
singularities with exponents $-1,0,1,2$. At infinity the exponents shift by 
$+1$ and are of the form
\[ 2+\alpha,\;\; 2-\alpha,\;\; 2+i \alpha,\;\;2-i \alpha \]
i.e. with $\sigma=4$ and $\Delta=-8$.

We remark that the exponents at infinity can only be all real in the very special case that these all coincide.\\

A precise statement relating these operators is the following:\\

\centerline{\Large \bf \em Proposition}
\vskip 8pt
\centerline{\em  $y(x)$ is a solution of $P_2(a,d,e,A)$}
\[ \Longleftrightarrow \]
\centerline{\em $(1+dx/4).y_0(x)$ is a solution of}
\vskip10pt
\centerline{$P_0(a,e-d/4,(4a-d)d/16,A-1)$}
 \[ \Longleftrightarrow \]
\centerline{\em $(1+(4a-d)x/4)^{-1}y(x)$ is a solution of}
\vskip 10pt
\centerline{ $P_4(a,e+(4a-d)/4,(4a-d)d/16), A-3)$}

\vskip 10pt

It is a remarkable fact that the local monodromies around the special points
$\neq 0$ and $\neq \infty$ are trivial. The monodromy around $0$ is MUM, so
has infinite order and thus the monodoromy group is the non-reductive group 
$\Z$, which clearly prevents these operators from being Picard-Fuchs operators 
of a family of algebraic varieties.\\

Something similar happens for the $\sigma=1$ and $\gamma=3$ components, given
by the operator families
\[
\begin{array}{rcl}
P_1(a,d,e,\alpha,\beta)&:=&\theta^4+x\left( a\theta^4-\frac{2d}{\Delta}\theta^3+\frac{\Delta-1}{\Delta}\theta^2+d\theta+e\right)+\\[1mm]
&&-\frac{d(d+a\Delta)}{\Delta^2}(\theta+\alpha)(\theta+\beta)(\theta+1-\alpha)(\theta+1-\beta)\\[1mm]
\end{array}
\]
where $\Delta=\alpha^2+\beta^2-\alpha-\beta$, and

\[
P_3(c,d,e,\alpha,\beta):=\theta^4 +x \left( \frac{(\Delta+8)c+(\Delta+9)d}{2(\Delta+6)}+2(c-d)\theta^3+c\theta^2+d\theta+e\right)+\]
\[+(2c-3d)\frac{(\Delta+4)c-(\Delta+3)d)}{2(\Delta+6)^2} x^2(\theta+\alpha)(\theta+\beta)(\theta+3+\alpha)(\theta+3-\beta),\]
with $\Delta=\alpha^2+\beta^2-3\alpha-3\beta$.\\

One verifies directly that\\

\centerline{\Large \bf \em Proposition}
\vskip 8pt
\centerline{\em $y(x)$ is a solution of $P_1(a,d,e,\alpha,\beta)$} 
\[ \Longleftrightarrow \]
\centerline{\em $(\Delta-(d+a\Delta)x)y_0(x)$ is a solution of}
\vskip10pt
\centerline{ $P_3(d (\Delta-1)/\Delta,d,e,\alpha+1,\beta+1)$}
\vskip 20pt

These operators are also remarkable. One of the roots of the polynomial
$1+ax+fx^2$ is in fact a regular point, the other is a conifold point, so 
the three non-trivial monodromies are precisely as those for the hypergeometric
operators. In fact, if we multiply the hypergeometric operator
\[\theta^4-Nx(\theta+\alpha)(\theta)\]\\
with the linear factor $1+Mx$, the result is the operator
\[P_1(M+N,-N\Delta,e,\alpha,\beta),\;\;\;(\Delta=\alpha^2+\beta^2-\alpha-\beta)\]
for the specific value 
\[ e=N\alpha\beta(1-\alpha)(1-\beta)\]
of the accessory parameter $e$. If $e$ takes on another value, the
operator is no longer equal to linear factor $\times$ hypergeometric.
\newpage

\section{\Large \bf \em The known Calabi-Yau operators of degree $2$}
\vskip 20pt
Below we will give an overview of the $69$ Calabi--Yau operators of degree two 
that are currently known to us. We also include some further remarkable
operators that are not strictly Calabi-Yau in the sense of $\cite{AZ}$.
Of course, this is just a small portion of the list of the still growing
list Calabi-Yau operators that was started in \cite{AESZ} and now contains
more then $500$ members. An update of this list in preparation \cite{ACvS}.
 
For operators of degree two there are cases with three and with four 
singular points. Most of the known operators can be related to hypergeometric 
operators or operators that are {\em convolutions}, that is, can be obtained 
via Hadamard product of operators of lower order. Recall that the {\em Hadamard product} 
of two power series
$\phi(x):=\sum_n a_n x^n$ and $\psi(x):=\sum_n b_n x^n$ is the series 
\[\phi \ast \psi (x) =\sum_n a_n b_n x^n\]
and by a classical theorem, if $\phi$ and $\psi$ satisfy a linear
differential equations, then so does $\phi \ast \psi$. The Hadamard product
corresponds to the multiplicative convolution of local systems, so
are under very good control, see \cite{DR}, \cite{DS}.\\

{\bf \em Frobenius basis:} For a fourth order operator $P$ with at MUM-point at
$0$, there is a canonical basis of solutions $y_0,y_1,y_2,y_3$ to $P y=0$ on a 
sufficiently small slit disc around the origin, called the {\em Frobenius basis}:
\[
\begin{array}{rcl}
y_0(x)&=&f_0(x)\\
&&\\
y_1(x)&=&\log(x) y_0(x)+f_1(x)\\
&&\\
y_2(x)&=&\frac{1}{2}\log(x)^2 y_0(x)+\log(x)y_1(x)+f_2(x)\\
&&\\
y_3(x)&=&\frac{1}{6}\log(t)^3y_0(x)+\frac{1}{2}\log(x)^2 y_1(x)+\log(x)y_2(x)+f_2(x)\\
\end{array}
\]
where $f_0(x) \in \Q[[x]]$, $f_i(x) \in t\Q[[t]], i=1,2,3$.\\

The Calabi-Yau operators from \cite{AZ} are characterised by three integrality
conditions:\\

1) integrality of the {\em solution} $y_0$: $$y_0(x) \in \Z[[x]].$$\\

2) integrality of the {\em q-coordinate} $q(x)$:\\

\[ q(x):=exp(y_1(x)/y_0(x))=x exp(f_1(x)/f_0(x))=x+\ldots \in \Z[[x]]\]

3) integrality of the {\em instanton numbers} $n_d$.\\

These can be computed in several different ways. 
The second logarithmic derivative of $y_2/y_0$ expressed in the $q$-coordinate
is called the {\em Yakawa-coupling} $K(q)$ of $P$:
\[K(q) = \left(q\frac{d}{dq}\right)^2 \left(\frac{y_2}{y_0} \right)\]
Now expand $K(q)$ as a Lambert-series
\[K(q)=1+\sum_{k=1}^{\infty} \frac{ k^3 n_k q^k}{1-q^k} .\]
The numbers $n_d$ are called the (normalised, $n_0=1$) {\em instanton numbers of $P$} and are required to to be integral. It is natural to allow small denominators in $f_0$, $q(x)$ and $n_d$, so require only $N$-integrality, for some denominator $N$.\\

\subsection{Description of the operators}
Many of the Calabi-Yau operators of degree 2 that we will describe involve 
special lower order operators of Calabi-Yau type. These are introduced and 
discussed in Appendix B.\\ 

{\large \bf \em Operators with three singular points}\\

{\em I. The 14 tilde operators}
There are $14$ exponents $(\alpha_1,\alpha_2,\alpha_3,\alpha_4)$ with
\[\alpha_1 \le \alpha_2 \le \alpha_3 \le \alpha_4,\;\;\alpha_1+\alpha_4=\alpha_2+\alpha_4=1,\;\] 
for which the hypergeometric operator, scaled by $N$,
\[\theta^4-Nx(\theta+\alpha_1)(\theta+\alpha_2)(\theta+\alpha_3)(\theta+\alpha_4)\]
is a Calabi-Yau operator, \cite{Alm1}. Corresponding to these, there are also 
$14$ hypergeometric {\em fifth order} Calabi-Yau operator
\[\theta^5-4Nx(\theta+\alpha_1)(\theta+\alpha_2)(\theta+\frac{1}{2})(\theta+\alpha_3)(\theta+\alpha_4)\]
with Riemann symbol
\[
\left\{
\begin{array}{ccc}
0&1/4N&\infty\\
\hline
0&0&\alpha_1\\
0&1&\alpha_2\\
0&3/2&1/2\\
0&2&\alpha_3\\
0&3&\alpha_4\\
\end{array}
\right\}
\]
These operators have a {\em Yifan Yang pull-back} to $14$ special fourth order 
operators, called the {\em tilde}-operators $\widetilde{1},\widetilde{2},\ldots, 
\widetilde{14}$. These operators replace the more complicated
{\em hat-operators}  $\hat{i},i=1,2,\ldots,14$ that appeared in the list \cite{AESZ} into which they can be transformed. However, operator $\widetilde{3}$, is not new, as it can be reduced to operator 2.33\\

{\em II. The 16 operators of type $H  \ast \mu(H)$}\\
 
There are four special hypergeometric second order operators called $A,B,C,D$.
The M\"obius transformation that interchanges the singularity with infinity 
leads to four {\em M\"obius transformed hypergeometric operators} 
\[ \mu(A):=e,\;\;\; \mu(B)=h,\;\;\; \mu(C)=i,\;\;\; \mu(D)=j ,\]
where $e,h,i,j$ refer to the names given in \cite{AZ}, see also Appendix B.\\

Thus we can form the $16$ Hadamard product $A \ast \mu(B)$, etcetera and obtain
$16$ fourth order operators of degree two with three singular points.\\
But there are two surprises: the operator $C \ast \mu(A)$ can be reduced to the 
hypergeometric operator $AESZ \#3=1.3$ and the operator $C \ast \mu(B)$ has vanishing 
instanton numbers and in fact is the third symmetric power of the hypergeometric 
second order operator 
\[{\theta}^{2}-12\, x\left( 12\,\theta+7 \right)  \left( 12\,\theta+1 \right)\]
So from this we obtain only $14$ Calabi-Yau operators of degree two that are 'new'.\\

{\em III. The four Hadamard products $I*\mu(H')$}\\ 

There are also four special hypergeometric third order operators
$A',B',C',D'$, analogous to second order operators $A,B,C,D$. Interchanging
the singularity with the point at infinity, we obtain third order operators
\[\mu(A')=\beta,\;\;\;\mu(B')=\iota,\;\;\;\mu(C')=\theta,\;\;\;\mu(D')=\kappa\] 
of degree two with three singular points. Taking Hadamard product with the 
central binomial coefficient produces four fourth order operators of degree 
two. One of these, $I \ast \beta$, can be transformd by a quadratic 
transformation to an operator $(I \ast \beta)^*$
\[\theta^4- 16(4\theta+1)(4\theta+3)(8\theta^2+8\theta+3)+4096x^2(4\theta+1)(4\theta+3)(4\theta+5)(4\theta+7),\]
which in turn can be transformed to the hypergeometric operator $AESZ \#3$. 
The operator $I \ast \theta$ can be transformed in $2.17$, so only two of these give rise to new Calabi-Yau operators.\\

{\bf IV. \em Four sporadic operators}\\ 

The operator $2.63=\#AESZ 84$ can not be constructed in this way. Furthermore, 
there are three operators of Bogner-type that were obtained from transformations 
of the higher degree operators $AESZ\#245, AESZ\#406$ and $AESZ \#255$. As these 
operators have unusual properties, they probably should better be understood in 
terms of the higher degree versions. We remark that the local monodromies 
around the two singularities $\neq 0$ are of finite order.\\

{\Large \bf \em Operators with four singular points}\\

{\em I. The $24$ operators of type $H \ast BZB$}
\vskip 10pt
There are six special second order operators of degree two with four singular 
fibres and Riemann symbol of the form
\[
\left\{  \begin{array}{cccc}0 &*&*&\infty\\
\hline
0&0&0&1\\
0&0&0&1\\
\end{array}\right\}.
\] 
So we are dealing here with special Heun operators. They appear in the work 
of Beukers \cite{Beukers1}, and Zagier \cite{Zagier}  and are related to the 
rational elliptic surfaces described by Beauville 
\cite{Beauville}, so we call them the BZB-operators. They were denoted
$a,b,c,d,f,g$ in \cite{AZ}.\\

Taking the Hadamard product with the hypergeometric operators $A,B,C,D$ produces 
$24$ fourth order operators of degree two with four singular points. \\

{ \em II. The six Hadamard products $I \ast BZB'$}\\

The BZB-operators have also analogs as third order operators that we will 
call the BZB'-operators, that were called $\alpha,\gamma,\delta,\epsilon,\zeta,\eta$ in \cite{AZ}. 
Taking $I \ast$ with these operators produce six new fourth order operators of degree two with four
singularities.\\

{\em III. Six Sporadic operators with four singularities}\\

{\em III.a. Three operators $I \ast$ sporadic third order}: There are three sporadic {\em third order} operators of degree two, described
in Appendix B. These lead to fourth order operators of degree two by taking 
$I \ast$ with them.\\ 

{\em III.b: Three original operators}: There are three operators that apparently can not be constructed starting from simpler operators of lower order; we call them original operators.\\

{ \bf \em The count}\\

So in total we have described
\[(14-1)+(16-2)+(4-2)+4=33\]
Calabi-Yau operators of degree two with three singular points and
\[24+6+3+3=36\]
with four singular points, making up a total of $69$. Because these all
have different instanton numbers, these are essentially different and there
are no algebraic transformations mapping one to any of the other operators.\\

{\bf \em Miscellaneous operators}\\

There are two further types of operators of degree two that deserve to be mentioned at this place, as they challenge the precise definition of Calabi-Yau operator, as formulated in \cite{AZ}. \\

{\em Reducible operators:}  A large scale Zagier-type search for degree two operators with integral solutions has been perfomed by {\sc Pavel Metelitsyn}. Apart from the operators equivalent to the ones mentioned above, there are
some interesting reducible fourth order operators. With regards to the q-coordinate  and instanton numbers these examples behave very much in the same way as the irreducible ones. An example is the operator
\[
{\theta}^{4}-12x\left( 6\theta+1 \right)\left( 6\theta+5 \right)\left( 2{\theta}^{2}+2\theta+1 \right) + 144 x^2 \left( 6\theta+1\right)\left( 6\theta+5 \right)  \left(6\theta+7\right)\left(6\theta+11\right) .\]
One finds:
\[\phi(x)=1 + 60 x + 13860 x^2  + 4084080 x^3  + 1338557220 x^4  + 465817912560 x^5+\ldots\]
\[q(x)= x + 312 x^2  + 107604 x^3  + 39073568 x^4  + 14645965026 x^5+\ldots\]  
\[n_1=-192,\;\;n_2=4182,\;\;n_3=-229568,\;\;n_4= 19136058,\;\; n_5=-2006581440\]

But note that the above series $\phi(x)$ is just the solution of the second order operator $D$, and the above operator is just the compositional square of of $D$. (Doing this with the hypergeometric second order operators $A,B,C$ produce
periodic instanton numbers.)\\

{\em A Strange Operator:} {\sc M. Bogner} \cite{Bogner1} also found an example of an operator of degree two with an integer solution, integral $q$-coordinate, integral $K(q)$-series, but for which the instanton numbers are badly non-integral. We start from the third order operator
\[
L:=\theta^{3}-4 x\left(2\theta+1 \right)\left(5\theta^{2}+5 \theta+2\right) 
+48 x^2 \left(3\theta+2 \right)\left(3 \theta+4 \right)\left(\theta +1 \right) 
\]
with integral solution
\[ 1 + 8 x + 96 x^2  + 1280 x^3  + 17440 x^4  + 231168 x^5+\ldots\]

The operator $I \ast L$ is the fourth order operator

\[\theta^4-8 x (2\theta + 1)^2(5\theta^2 + 5\theta +2) + 192 x^2(2\theta+ 1)(3\theta + 2)(3\theta + 4)(2\theta + 3) \]
We find the following nice integral series:
\[y_0(x)=1 + 16 x + 576 x^2  + 25600 x^3  + 1220800 x^4  + 58254336 x^5+\ldots\]
\[ q(x):= x + 40 x^2  + 1984 x^3  + 106496 x^4  + 5863168 x^5  + \ldots\]
\[K(q):= 1 + 8 q^2 - 5632 q^3  - 456064 q^4  - 17708032 q^5+ \ldots \]
However, if we expand $K(q)$ as a Lambert-series we do not get integral instanton numbers: 
\[n_1=8,n_2=-1,n_3=-{\frac {1880}{9}},n_4=-7126,n_5=-{\frac {3541608}{25}},\ldots,\]
\[n_{17}=\frac{2432475693294458880448632}{289},\ldots\]
Experimentally, the denominator $p^2$ appears in $n_p$ for
\[ p = 3,\;\; 5,\;\; 7,\;\; 11,\;\; 13\;\; 17\;\; 19,\ldots\]

This is rather puzzling. Conjecturally, already the fact that $\phi(x)$ is an integral series implies that the operator is a factor of a Picard-Fuchs equation of a family of varieties defined over $\Q$ and indeed our operator is of geometrical origin. It was conjectured in \cite{AZ} that the integrality of the $q$-series implies the integrality of the instanton numbers, which in this example is not the case. The integrality of solution and mirror map clearly indicate that we have a rank four Calabi-Yau motive and one would expect the general arguments for the integrality of \cite{Vologodsky} to be applicable, but apparently they are not. There might exist a different scaling of the coordinate that repairs this defect, but up to now we have been unable to find it.\\

{\em \bf Acknowledgement} This work for this paper was done in 2016 and we
had hoped to find many new degree two operators. That did not happen and we
never found the time to write the paper. We intended it to be a possible
reference for this class of operators. The untimely death of Gert left the
paper in a unfinished state and it is only now I found the time to complete 
it. We thank M. Bogner, J. Hofmann and K. Samol for help with computations, done
years ago, and W. Zudilin and O. Gorodetsky for useful comments, corrections 
and showing interest in this paper.\\ 

\section*{\Large \bf \em Appendix A}
The following table lists the operators as they appear in the 
electronic database accessable at {\tt https://cydb.mathematik.uni-mainz.de/}.
(These are ordered by degree, so that the $n$th operator of degree $d$ gets the
code $d.n$ in that list.) 
The second column lists the AESZ-numbers as appearing in \cite{AESZ}. The last 
three list the first three instanton numbers (normalised, $n_0=1$). As these
do not change under transformations of the coordinate, the pair of numbers 
$|n_1|,|n_3|$ provide a simple 'fingerprint/superseeker' to identify an
operator. In Appendix C one finds some monodromy and modular form information 
for these operators.\\ 

\[
\begin{array}{|c|c|c|c||r|r|r|}
\hline
Number&AESZ      &Source  &|\Sigma|&n_1&n_2&n_3\\
\hline
2.1    &\#45      &A\ast a &4  &12&163 &3204\\
\hline
2.2    &\#15      &B \ast a&4  &21&480 &15894\\
\hline
2.3    &\#68      &C \ast a&4  &52&2814 &220220\\ 
\hline
2.4    &\#62      &D \ast a&4  &372&136182 &71562236\\     
\hline
2.5    &\#25      &A \ast b&4  &20&277 &8220\\      
\hline
2.6    &\#24      &B \ast b&4  &36&837 &41421\\
\hline
2.7    &\#51      &C \ast b&4  &92&5052 &585396\\
\hline
2.8    &\#63      &D \ast b&4  &684&253314 &195638820\\
\hline
2.9    &\#58      &A \ast c&4  &16&142&11056/3 \\
\hline
2.10   &\#70      &B \ast c&4  &27&432&18089\\  
\hline
2.11   &\#69      &C \ast c&4  &64& 2616& 246848\\
\hline
2.12   &\#64      &D \ast c&4  &432&130842&78259376\\
\hline
2.13   &\#36      &A \ast d&4  &16& 42& 1232\\
\hline
2.14   &\#48      &B \ast d&4  &24& 291/2& 5832\\
\hline
2.15   &\#38      &C \ast d&4  &48& 998& 73328 \\
\hline
2.16   &\#65      &D \ast d&4  &240& 57102& 19105840 \\      
\hline
2.17   &\#111,\sim\#40 &A \ast e&3  &32& -96& 1440\\
\hline
2.18   &\#110     &B \ast e&3  &36& -144& 8076\\
\hline
2.xx19 & \sim \#3 &C \ast e&3  & 32& 608& 26016 \\
\hline
2.19   &\#112     &D \ast e&3  &-288& 162504& -96055968 \\
\hline
2.20   &\#133        &A \ast f&4  &12& -42& -3284/3 \\
\hline
2.21   &\#134        &B \ast f&4  &18& -207/2& -5177  \\
\hline
2.22  &\#135         &C \ast f&4  &36& -477& -206716/3 \\
\hline
2.23  &\#136         &D \ast f&4  &180& -15615& -21847076   \\
\hline
2.24  &\#137        &A \ast g&4  &20& 2& 1684/3    \\
\hline
2.25  &\#138        &B \ast g&4  &27& 189/4& 2618  \\ 
\hline
2.26  &\#139        &C \ast g&4  &44& 607& 22500   \\  
\hline
2.27  &\#140        &D \ast g&4  &108& 54135& -4945756 \\
\hline
2.28  &\#141        &A \ast h&3  & 48& -438& 2864    \\
\hline
2.29  &\#142        &B \ast h&3  & 45& -3465/4& 27735    \\
\hline
2.xx30& Sym^3         &C \ast h&3  &  0& 0& 0  \\
\hline
2.30&\#143          &D \ast h&3  &  -1008&499086& -607849200 \\
\hline
\end{array}
\]
\[
\begin{array}{|c|c|c|c||r|r|r|}
\hline
Number&AESZ      &Source  &   &n_1&n_2&n_3\\
\hline
2.31&  7^{**}      &A \ast i&3  &96& -3560& -12064\\
\hline
2.32&            &B \ast i&3  & 60& -7635& 307860 \\
\hline
2.33& 6^*          &C \ast i&3  &-160& -6920& -539680\\
\hline
2.34&            &D \ast i&3  &-3936& 3550992& -10892932064\\
\hline
2.35& 9^{**}, \sim \#67&A \ast j&3  &  480& -226968& -16034720\\
\hline
2.36&            &B \ast j&3  &-36& -486279& 128217204\\
\hline
2.37&            &C \ast j&3  & -2592& -307800& 81451104\\
\hline
2.38& \#61       &D \ast j&3  &  -41184& 251271360& -5124430612320\\
\hline
2.39& \sim \widehat{1}& \widetilde{1}&3& -2450& -1825075/2& -623291900\\
\hline
2.40&\sim \widehat{2}& \widetilde{2}&3&-791200& -41486886600&-4288711075194400 \\\hline
2.xx41& \sim \widehat{3} & \widetilde{3}&3&-160& -6920& -539680\\
\hline
2.41&\sim \widehat{4}& \widetilde{4}&3&-522& -105291/2& -9879192\\
\hline
2.42&\sim \widehat{5}& \widetilde{5}&3& -288& -19260& -2339616\\
\hline
2.43&\sim \widehat{6}& \widetilde{6}&3&-736& -104512& -26911072\\
\hline
2.44&\sim \widehat{7}& \widetilde{7}&3&-57760& -354010600& -3869123234080\\
\hline
2.45&\sim \widehat{8}& \widetilde{8}&3&-10080& -11338740& -24400330080\\
\hline
2.46&\sim \widehat{9}& \widetilde{9}&3& -2710944& -717640301160& -302270555492914464\\
\hline
2.47&\sim \widehat{10}& \widetilde{10}&3&-3488& -1406056& -1142687008 \\
\hline
2.48&\sim \widehat{11}& \widetilde{11}&3&-1344& -278040& -109320512 \\
\hline
2.49&\sim \widehat{12}& \widetilde{12}&3&-264006& -52511160& -230398034080\\
\hline
2.50& \sim \widehat{13}& \widetilde{13}&3&-201888& -1567499400& -40177844666400\\
\hline
2.51&\sim \widehat{14}& \widetilde{14}&3&-5472& -4476528& -6444589536\\
\hline
2.52 &\#16&I \ast \alpha&4& 4& 20& 644/3  \\
\hline
2.53 &\#29&I \ast \gamma&4& 14& 303/2& 10424/3 \\
\hline
2.54& \#41&I \ast \delta&4& 2& -7& -104 \\
\hline
2.55 &\#42&I \ast \epsilon&4&8& 63& 1000\\
\hline
2.56 &\#185&I \ast \zeta&4&6& 93/2& 608 \\
\hline
2.57 &\#184&I \ast \eta&4& 2& 4& -8 \\
\hline
2.xx58&\sim \#3&I \ast \beta&3& 0& 4& 0 \\
\hline
2.yy58&\sim 2.xx58&(I\ast \beta)^*&3&32& 608& 26016\\
\hline
2.58 &4*  &I\ast \iota&3 & -6& -6& -104 \\
\hline
2.xx59&\sim 2.17&I \ast \theta&3& -32& -88& -1440 \\
\hline
2.59 & 13^*,\sim 47 &I \ast \kappa&3 &-384& -1356& -164736\\
\hline
2.60 &\#18  &I \ast Spor1   &4  & 4& 39& 364  \\
\hline
2.61 &\#26  &I \ast Spor2  &4   & 10& 191/2& 1724  \\
\hline
2.62 &\#28  &Original   &4   & 5& 28& 312  \\
\hline
2.63 &\#84  &Original   &3   & -4& -11& -44 \\
\hline
2.64 &\#182 &Original   &4   & 1& 7/4& 7  \\
\hline
2.65 &\#183 &I \ast Spor3  &4   & 4& 7& 556/9  \\
\hline
2.66&reducible!& D  \circ D &3   & -192& 4182& -229568\\
\hline
2.67&  \sim \#245 &Bogner1   &3   & 6& -69/2& 170  \\
\hline
2.68&  \#406 &Bogner2       &3   & -12& -186& -1668  \\
\hline
2.69& \#205 &Original       &4  &   1& 7/4& 5  \\ 
\hline 
2.70&\sim \#255&Bogner3     &3  &  20& 290& 28820/3\\
\hline
\end{array}
\]

\section*{\Large \bf \em Appendix B}

In this appendix we collect some information on operators of lower order 
that are relevant  to the construction of the order four Calabi-Yau operators 
of degree two.

\subsection*{\Large \bf \em First order operator}

There is a single {\em first order operator} that we call $I$:

\[ \theta-4x(\theta+1/2)\]
The series
\[\sum_{n=0}^{\infty} {2 n \choose n} x^n\]
is the unique holomorphic solution of $I$, and in fact coincides for
$x \le 1/4$ with the algebraic function
\[\frac{1}{\sqrt{1-4x}} .\]

\subsection*{ \Large \bf \em Second order operators}
\vskip 10pt
\subsubsection*{\large \bf \em The four hypergeometric cases}

There are four very remarkable hypergeometric second order operators 
of the form
\[ \theta^2 -N x (\theta+\alpha)(\theta+\beta),\;\;\;\alpha+\beta=1,\;\;\;\alpha=1/2,\;1/3,\;1/4,\;1/6\,\]
with Riemann symbol of the form
\[
\left\{
\begin{array}{ccc}
0&1/N&\infty\\
\hline
0&0&\alpha\\
0&0&\beta\\
\end{array}
\right\} .
\]
These four operators appear in many different areas of mathematics, like 
the {\em alternative theories of Ramanujan} and are described at many places, 
see e.g. \cite{Cooper}.\\

\[
\begin{array}{|c|l|l|c|}
\hline
\textup{Name}&\textup{Operator}&a_0,a_1,a_2,a_3,\ldots &a_n\\
\hline
A &\theta^2-16 x (\theta+\frac{1}{2})^2&1, 4, 36, 400,\ldots &\frac{(2n!)^2}{n!^4}\\[1mm]
\hline
B&\theta^2-27 x (\theta+\frac{1}{3})(\theta+\frac{2}{3})& 1, 6, 90, 1680, \ldots &\frac{(3n)!}{n!^3}\\[1mm]
\hline
C&\theta^2-64 x (\theta+\frac{1}{4})(\theta+\frac{3}{4}) & 1, 12, 420, 18480,\ldots  &\frac{(4n)!}{(2n!) (n!)^2}\\[1mm]
\hline
D&\theta^2-432 x (\theta+\frac{1}{6})(\theta+\frac{5}{6})& 1, 60, 13860, 4084080, \ldots &\frac{(6n)!}{(3n!)(2n!)(n!)}\\[1mm]
\hline
\end{array}
\]

The above list of operators is in close correspondence with the
six extremal elliptic surfaces with three singular fibres, \cite{Schmickler}, \cite{MP}.

\[
\begin{array}{|c|c|c|c|c|c|c|c|}
\hline
\textup{fibres}&\textup{MW}&A&B&(BA)^{-1}&\textup{Monodromy}&\textup{H/MP}&\\
\hline
I_2 \;\;I_2\;\;I_2^*&(\Z/2)^2&\left (
\begin{array}{rr}
1&2\\
0&1\\
\end{array} \right )
& \left (
\begin{array}{rr}
1&0\\
-2&1\\
\end{array} \right )
&
\left (
\begin{array}{rr}
1&-2\\
2& -3\\
\end{array} \right )
&\Gamma(2)  & 12 /X_{222}&A\\
\hline
I_1\;\;I_1\;\;I_4^*&\Z/2
&
\left (
\begin{array}{rr}
1&1\\
0&1\\
\end{array} \right )

& 
\left (
\begin{array}{rr}
-1&1\\
-4&3\\
\end{array} \right )

& 
\left (
\begin{array}{rr}
-1&0\\
 4&-1\\
\end{array} \right )

&\overline{\Gamma_0(4)} &11/X_{114}&A\\
\hline
\hline
I_1\;\;I_4\;\;I_1^*&\Z/4&\left (
\begin{array}{rr}
1&1\\
0&1\\
\end{array} \right ) &\left (
\begin{array}{rr}  
1&0\\
-4&1\\
\end{array} \right ) &\left (
\begin{array}{rr}
-3&-1\\
 4&1\\
\end{array} \right ) &\Gamma_0(4) &11/X_{141}&A\\
\hline
I_1 \;\;I_3\;\;IV^*&\Z/3&
\left (
\begin{array}{rr}
1&1\\
0&1\\
\end{array} \right )
& \left (
\begin{array}{rr}
1&0\\
-3&1\\
\end{array} \right )
&
\left (
\begin{array}{rr}
1&-1\\
3& -2\\
\end{array} \right )&\Gamma_0(3)   &10/X_{134}&B\\
\hline
I_1 \;I_2\;\; III^*&\Z/2& 

\left (
\begin{array}{rr}
1&1\\
0&1\\
\end{array} \right )

& \left (
\begin{array}{rr}
1&0\\
-2&1\\
\end{array} \right )&

\left (
\begin{array}{rr}
1&-1\\
2&-1\\
\end{array} \right )
&\Gamma_0(2)&9/X_{123}&C\\
\hline
I_1 \;I_1\;\; II^*&\Z/1&

\left (
\begin{array}{rr}
1&1\\
0&1\\
\end{array} \right )

& \left (
\begin{array}{rr}
1&0\\
-1&1\\
\end{array} \right )&

\left (
\begin{array}{rr}
1&-1\\
1& 0\\
\end{array} \right )
 &\Gamma_0(1)& 8 /X_{112}&D\\
\hline
\end{array}
\]
Here $MW$ stands for the Mordell-Weil group of sections of the fibration.
The matrices $A$, $B$ and $(BA)^{-1}$ are monodromy matrices for
loops around $0$, $c\neq 0, \infty$ and $\infty$ respectively. In the
monodromy column, the group generated by $A$ and $B$ is identified. 
The column $H/MP$ denotes the notation used in \cite{Schmickler} and \cite{MP}. 
The last column denotes the type of the Picard-Fuchs operator. We see that
there are three surfaces that belong to operator $A$, but the surface
$I_1\;\;I_4\;\;I_1^*$ shows it is a member of a very regular series of
four surfaces, which we will say to belong to the operators $A,B,C,D$. 

The surfaces belonging to the operators $A,B,C,D$ have a clear
relation to the geometry of the {\em simple elliptic singularities}.
The fibres over $\infty$ correspond to configurations of rational
curves intersecting in the pattern of the affine Dynkin diagrams.
The corresponding simple elliptic singularities have a single modulus.
These families of elliptic curves lead to elliptic surfaces that cover
the corresponding extremal elliptic surfaces.
\[ 
\begin{array}{|c|c|c|c|}
\hline
\textup{Case}&\textup{Dynkin}&\textup{Equation}&\textup{Cover}\\[1mm]
\hline
A&\widetilde{D_5}&x^2+y^2=suv,\;\;u^2+v^2=sxy&t=s^2\\
B&\widetilde{E_6}&x^3+y^3+z^3=sxyz&t=s^3\\
C&\widetilde{E_7}&x^2+y^4+z^4=sxyz&t=s^4\\
D&\widetilde{E_8}&x^2+y^3+z^6=sxyz&t=s^6\\
\hline
\end{array}
\]
The $A$ case is slightly different, as the corresponding singularity
is not a hypersurface, but a complete intersection of two quadrics.
The last three cases correspond to the {\em euclidian triples}
$(p,q,r)$ of integers $\ge 2$ with the property that
\[ \frac{1}{p}+\frac{1}{q}+\frac{1}{r}=1 .\]

\subsubsection*{\bf \em The transformed hypergeometric cases}
 
The M\"obius-transformation 
\[ x \mapsto - \frac{x}{1-Nx}\]
preserves $0$ and interchanges the points $1/N$ and $\infty$ of the
Riemann sphere. If $y(x)$ solves one of the hypergeometric equations $A,B,C,D$, then
\[ Y(x):=\frac{1}{1-Nx} y(-\frac{x}{1-Nx})\]
satisfies a second order equation of the form
\[ \theta^2-x(a\theta^2+a\theta+b)+cx^2(\theta+1)^2\]
that we call {\em transformed hypergeometric equation} that we
denote by $\mu(A),\mu(B),\mu(C),\mu(D)$. These operators have 
Riemann symbol of the form
\[
\left\{
\begin{array}{ccc}
0&1/N&\infty\\
\hline
0&-\alpha&1\\
0&-\beta&1\\
\end{array}
\right\}
\]

\[
\begin{array}{|c||c|c|c||l|c|}
\hline
\textup{Name}&a&b&c&a_0,a_1,a_2,\ldots &a_n\\[1mm]
\hline
\mu(A)=e&32&12&16^2&1,12,164,\ldots&16^n \sum_{k=0}^n (-1)^k{-1/2 \choose k}^2{-1/2 \choose n-k}\\[1mm]
\hline
\mu(B)=h&54&21&27^2&1,21,495,\ldots&27^n \sum_{k=0}^n (-1)^k {-1/3 \choose k}^2{-2/3 \choose n-k}\\[1mm]
\hline
\mu(C)=i&128&52&64^2&1,52,2980,\ldots&64^n \sum_{k=0}^n (-1)^k {-1/4 \choose k}^2{-3/4 \choose n-k}\\[1mm]
\hline
\mu(D)=j&864&372&432^2&1, 372, 148644,\ldots&432^n \sum_{k=0}^n (-1)^k {-1/6 \choose k}^2{-5/6 \choose n-k}\\[1mm]
\hline
\end{array}
\]

\subsubsection*{\large \bf \em The six Beukers-Zagier-Beauville operators}
\vskip 10pt
There are six special second order operators of degree two with four
singular points of the form
\[ \theta^2-t(a\theta^2+a\theta+b)+c(\theta+1)^2\]
and hence have a Riemann-symbol of the form
\[
\left\{
\begin{array}{cccc}
0&*&*&\infty\\
\hline
0&0&0&1\\
0&0&0&1\\
\end{array}
\right\}
\]

\[
\begin{array}{|r|r||r|r|l|l|c|}
\hline
\textup{AZ-name}&\textup{BZ-name}&a&b&c&a_0,a_1,a_2,a_3, \ldots&a_n \\[1mm]
\hline
\hline
a&A&7&2&-8&1, 2, 10, 56, \ldots  &\sum_{k=0}^n {n \choose k}^3\\[3mm]
\hline
c&C&10&3&9&1, 3, 15, 93, \ldots &\sum_{k=0}^n {n \choose k}^2 {2k \choose k}\\[3mm]
\hline
g&F&17&6&72&1, 6, 42, 312, \ldots &\sum_{j,k=0}^n (-1)^j 8^{n-j}{n \choose j} {j \choose k}^3\\[3mm]
\hline
\hline
d&E&12&4&32&1, 4, 20, 112,\ldots &\sum_{k=0}^n {n \choose k} {2k \choose k}{2n-2k \choose n-k}\\[3mm]
\hline
f&B&9&3&27&1, 3, 9, 21, \ldots &\sum_{k=0}n (-1)^k3^{n-2k}{n \choose 3k }\frac{(3k)!}{k!^3}\\[3mm]
\hline
b&D&11&3&-1&1, 3, 19, 147,\ldots &\sum_{k=0}^n {n \choose k}^2{n+k \choose k}\\[3mm]
\hline
\end{array}
\]
We refer to \cite{Gorodetsky} for a representation of these sequences
as constant term of a  Laurent polynomial.\\

Unfortunately there is no natural naming for these operators. In the first column of the
table we used the ``AZ-names'' used in \cite{AZ}. In  \cite{Beukers2} and \cite{Zagier} 
the  same operators are named $A,B,\ldots,F$ (Not to be confused with the hypergeometric
operators $A,B,C,D$ mentioned above!); these names appear in the second column as
``BZ-names''. These operators also appear as Picard-Fuchs operators for rational 
elliptic surfaces, namely the six semi-stable families of elliptic curves with four 
exceptional fibres,  that have been studied by {\sc Beauville}, \cite{Beauville}. 
The modular level gives a more or less natural way order the operators.\\

\[
\begin{array}{|c|c|c|c|c|}
\hline
\textup{fibres}&j&Pencil&Group&Operator\\[1mm]
\hline
\hline
I_5\;I_5\;I_1\;I_1&2^{14}31^3/5^3&x(x-z)(y-z)+tzy(x-y)&\Gamma_1(5)& b  \\[1mm] 
\hline
I_6\;I_3\;I_2\;I_1&2^2.73^3/3^4 &(x+y)(y+z)(z+x)+txyz&\Gamma_0(6)&a,c,g\\[1mm]
\hline
I_8\;I_2\;I_1\;I_1&1728& (x+y)(xy-z)+txyz &\Gamma_0(8) \cap \Gamma_1(4)& d \\[1mm]
I_4\;I_4\;I_2\;I_2& &x(x^2+z^2+2yz)+tz(x^2-y^2)&\Gamma_1(4) \cap \Gamma(2)&   \\[1mm]
\hline
I_9\;I_1\;I_1\;I_1&0&x^2y+y^2z+z^2x+txyz&\Gamma_0(9) \cap \Gamma_1(3)&  f\\[1mm]
I_3\;I_3\;I_3\;I_3& &x^3+y^3+txzy& \Gamma(3) &        \\[1mm]
\hline
\end{array}
\]
The coefficients of the holomorphic solution of the operator $b$ 
are the Ap\'ery-numbers for the irrationality of $\pi^2$.
It belongs to the elliptic surface with fibres $I_5,I_5,I_1,I_1$
over the modular curve $X_1(5)$, as was first noticed by Beukers in
\cite{Beukers1}.\\
The operators $a,c,g$ differ by a M\"obius-transformations.
The base of this family can be identified with the modular curve for 
the congruence subgroup $\Gamma_0(6)$. Over that curve we have the 
semi-stable elliptic surface with fibres of type $I_6,I_3,I_2,I_1$. 
The operator $d$ belongs to the two isogenous elliptic surfaces
with fibres $I_8,I_2,I_1,I_1$ and $I_4,I_4,I_2,I_2$. Operator $f$ belongs 
to the two isogenous elliptic surfaces
with fibres $I_9,I_1,I_1,I_1$ and  $I_3,I_3,I_3,I_3$ which is nothing
but the Hesse pencil.\\

\subsection*{\Large \bf \em Third order operators}

\subsubsection*{\large \bf \em Hypergeometric operators}
Closely related to the operators $A,B,C,D$ there are also four very remarkable hypergeometric third order operators 
\[A'=I \ast A,\;\;\;B'=I \ast B,\;\;\;C'=I \ast C,\;\;\;D'=I \ast D\] of the form
\[ \theta^3 -N x (\theta+\alpha)(\theta+1/2)(\theta+\beta),\;\;\;\alpha+\beta=1,\;\;\;\alpha=1/2,\;1/3,\;1/4,\;1/6\,\]
with Riemann symbol
\[
\left\{
\begin{array}{ccc}
0&1/N&\infty\\
\hline
0&0&\alpha\\
0&1/2&1/2\\
0&0&\beta\\
\end{array}
\right\} .
\]

Also important are the {\em symmetric squares} of the operators $A,B,C,D$.
These are of the form
\[ \theta^3-x(2\theta+1)(a\theta^2+a\theta+b)+a^2 x^2(\theta+\alpha)(\theta+1)(\theta+\beta),\;\;\alpha+\beta=2\]
and have Riemann-symbol of the form
\[
\left\{
\begin{array}{ccc}
0&1/a&\infty\\
\hline
0&0&1-\alpha\\
0&0&1\\
0&0&1+\alpha\\
\end{array}
\right\}
\]

\[
\begin{array}{|c||c|c|c|}
\hline
\textup{Name}&a&b&\alpha\\
\hline
sym^2(A)&16&8&1-2\frac{1}{2}=0\\
\hline
sym^2(B)&27&12&1-2\frac{1}{3}=1/3\\
\hline
sym^2(C)&64&24&1-2\frac{1}{4}=1/2\\
\hline
sym^2(D)&432&120&1-2\frac{1}{6}=2/3\\
\hline
\end{array}
\]

\subsubsection*{\large \bf \em The four transformed symmetric squares}
\vskip 10pt
More important for us are the {\em transformed symmetric squares} of the operators 
$A,B,C,D$. These are of the form
\[ \theta^3-x(2\theta+1)(a\theta^2+a\theta+b)+a^2x^2(\theta+1)^3\]
and have Riemann-symbol of the form
\[
\left\{
\begin{array}{ccc}
0&1/a&\infty\\
\hline
0&-\alpha&1\\
0&0&1\\
0&\alpha&1\\
\end{array}
\right\}
\]

\[
\begin{array}{|c||c|c|c||c|c|}
\hline
\textup{Name}&a&b&\alpha&a_0,a_1,a_2,\ldots&a_n\\
\hline
\mu(sym^2(A))=\beta&16&8&1&1, 8, 88, 1088,  &16^n\sum_k {-1/2\choose k}^2{-1/2 \choose n-k}^2\\
\hline
\mu(sym^2(B))=\iota&27&15&1/3&1, 15, 297, 6495&27^n\sum_k {-1/3\choose k}^2{-2/3 \choose n-k}^2\\
\hline
\mu(sym^2(C))=\theta&64&40&1/2&1, 40, 2008, 109120&64^n\sum_k {-1/4\choose k}^2{-4/3 \choose n-k}^2\\
\hline
\mu(sym^2(D))=\kappa&432&312&2/3&1, 312, 114264, 44196288&432^m\sum_k {-1/6\choose k}^2 {-5/6 \choose n-k}^2\\
\hline
\end{array}
\]
The names $\beta, \iota, \theta,\kappa$ stem from \cite{AZ}.

\subsubsection*{\large \bf \em The six BZB' -operators}
\vskip 10pt
Corresponding to the six BZB-operators, there are also 
six special third order operators of degree two with four
singular points of the form
\[ \theta^3-x(2 \theta+1)(a\theta^2+a\theta+b)+cx^2(\theta+1)^3 \]
and hence have a Riemann-symbol of the form
\[
\left\{
\begin{array}{cccc}
0&*&*&\infty\\
\hline
0&0&0&1\\
0&1/2&1/2&1\\
0&1&1&1
\end{array}
\right\}
\]

\[
\begin{array}{|r|r|r|l|l|c|}
\hline
\textup{Name}&a&b&c&a_0,a_1,a_2,a_3, \ldots&a_n \\
\hline
\delta&7&3&81&1, 3, 9, 3, -279&\sum_k (-1)^k 3^{n-3k} {n \choose 3k}{n+k \choose k}\frac{(3k)!}{k!^3}\\[2mm]
\hline
\alpha&10&4&64&1, 4, 28, 256,\ldots &\sum_k{n \choose k}^2{2k \choose k}{2n-2k \choose n-k}\\[2mm]
\hline
\gamma&17 & 5& 1&1, 5, 73, 1445, \ldots &\sum_k {n \choose k}^2{n+k \choose k}^2 \\[2mm]
\hline
\hline
\epsilon&12 & 4& 16&1, 4, 40, 544, \ldots &\sum_k {n \choose k}^2{2k \choose n}^2\\[2mm]
\hline
\eta&11&5&125&1, 5, 35, 275,\ldots & \sum_{k=0}^{\lfloor n/5 \rfloor} {n \choose k}^3\left({4n-5k-1\choose 3n}+{4n-5k \choose 3n} \right) \\[2mm]
\hline
\zeta&9& 3& -27&1, 3, 27, 309, \ldots &\sum_{k,l} {n \choose k}^2{n \choose l}{k \choose l}{k+l \choose n} \\[2mm]
\hline
\end{array}
\]
For a Laurent-representation  of these sequences we refer to \cite{Gorodetsky}.\\

The correspondence between the $BZB$-operators and $BZB'$-operators is determined 
by the following table:
{\Large
\[
\begin{array}{|c|c|c||c|c|c|}
\hline
a&c&g&d&f&b\\
\hline
\delta&\alpha&\gamma&\epsilon&\zeta&\eta\\
\hline
\end{array}
\]
}
There is a general formula expressing the relation between the
BZB and BZB' operators:

\[  Y_0(-x/Q(x))=Q(x)y_0(x)^2\]

which shows that the third order operators are twisted versions of the symmetric squares
of the corresponding second order operators. 
The factor $Q(x)$ is the discriminant of the second order operator with $y_0$ as solution.
\[
\begin{array}{|c|c|c|}
\hline
y_0&Y_0&Q\\
\hline
a&\delta&(1+x)(1-8x)\\
c&\alpha&(1-x)(1-9x)\\
g&\gamma&(1-8x)(1-9x)\\
\hline
d&\epsilon&(1-4x)(1-8x)\\
\hline
f&\zeta&1-9x+27x^2\\
\hline
b&\eta&1-11x-x^2\\
\hline
\end{array}
\]

\subsubsection*{\large \bf \em The three sporadic third order operators}
\vskip 10pt
We know of three very special third order operators that seem to be not 
directly to the analogous of the second order operators. These are the 
following ones:\\ 

{\bf Sporadic 1:}

\[\theta^3-2x(2\theta+1)(3\theta^2+3\theta+1)-4x^2(4\theta+3)(\theta+1)(4\theta+5)\]
\[\phi(x)=1+ 2x +18x^2 +164x^3+ 1810x^4+21252x^5+263844x^6+\ldots \]

\[
\left\{
\begin{array}{cccc}
0&1/16&-1/4&\infty\\
\hline
0&0   &0  &3/4\\
0&1/2 &1/2&1\\
0&1   &1  &5/4\\
\end{array}
\right\}
\]
\[a_n =\sum_{k=0}^n  {n  \choose k}^4\]
This example has an A-incarnation as K3-surface that is the intersection of four
hyperplane sections of type $(1,1)$ in $\P^3 \times \P^3$. If we take $I \ast$, we obtain a fourth order operator 2.60 with A-incarnation a CY threefold that appears as  intersection of two $(1,1)$ and one $(2,2)$ hypersurface in $\P^3 \times \P^3$, \cite{BvS}.\\

{\bf Sporadic 2:}
\[\theta^3-x(2\theta+1)(13\theta^2+13\theta+4)-3x^2(3\theta+2)(\theta+1)(3\theta+4)\]
\[\phi(x)=1+ 4 x+ 48 x^2+760 x^3+13840 x^4+ 273504 x^5+5703096 x^6+\ldots\]
\[
\left\{
\begin{array}{cccc}
0&1/27&-1&\infty\\
\hline
0&0   &0  &2/3\\
0&1/2 &1/2&1\\
0&1   &1  &4/3\\
\end{array}
\right\}
\]

\[a_n =\sum_{k=0}^n  {n  \choose k}^2 {n+k \choose n}{2k \choose n}\]

This example has an A-incarnation as K3-surface that is the intersection of six
hyperplane sections of the Grassmanian $G(2,6)$ in its Pl\"ucker embedding. 
If we take $I \ast$, we obtain a fourth order operator 2.61 with A-incarnation a 
CY threefold that appears as intersection of four linear and a quadratic hypersurface 
in Grassmanian $G(2,6)$ in its Pl\"ucker embedding, \cite{BCKvS}.\\

{\bf Sporadic 3:}
\[\theta^3-2x(2\theta+1)(7\theta^2+7\theta+3)+12x^2(4\theta+3)(\theta+1)(4\theta+5)\]
\[\phi(x)= 1+6x+ 54x^2+564x^3+6390x^4+76356x^5+948276x^6+\ldots\]

\[
\left\{
\begin{array}{cccc}
0&1/16&1/12&\infty\\
\hline
0&0   &0  &3/4\\
0&1/2 &1/2& 1\\
0&1   &1  &5/4\\
\end{array}
\right\}
\]

\[a_n =\sum_{k=0}^{\lfloor n/3 \rfloor} (-1)^k {n \choose k} {2k \choose k}{2(n-k) \choose n-k}\left({2n-3k-1 \choose n}+{2n-3k \choose n}\right)
\]
\\

These sporadic operators and sequences were also found by S. Cooper 
\cite{Cooper0}, where they are called $s_{10}, s_{7}$ and $s_{18}$ and where
the explicit description for $a_n$ in the last case also can be found.
We refer to \cite{Gorodetsky} for representations as constant term sequence
of a Laurent polynomial.\\

\subsubsection*{ \bf \em Miscellaneous third order operators}
\vskip 10pt
There is a remarkable infinite series of third order operators with an integral
solution.\\

\centerline{\Large \bf \em Theorem}
\vskip 10pt 
{\em The differential operator 
\[P:=\theta^3-t(2\theta + 1)(\theta^2  + \theta - 1/2d^2  + 1/2k^2  + 1/2)+x^2 (\theta + 1 - d)(\theta + 1)(\theta + 1 + d)\]
has 
\[(1-x)^k \cdot  _2F_1\left(\frac{1-d+k}{2},\frac{1+d+k}{2},1, x\right)^2\]
as solution.}

{\bf \em Proof:} The Riemann symbol of the operator $P$ is

\[
\left\{ \begin{array}{ccc}
0&1&\infty\\
\hline
0&-k&(1-d)\\
0&0& 1\\.
0&k&1+d\\
\end{array}\right\} 
\]

The operator is the symmetric square of a second order operator
with Rieman symbol

\[
\left\{ \begin{array}{ccc}
0&1&\infty\\
\hline
0&-k/2&(1-d)/2\\
0&k/2&(1+d)/2 \\
\end{array}\right\} .
\]
Multiplication by a factor $(1-x)^{k/2}$ gives a hypergeometric 
operator with Riemann symbol
\[
\left\{ \begin{array}{ccc}
0&1&\infty\\
\hline
0&0&(1-d+k)/2\\
0&k&(1+d+k)/2 \\
\end{array}\right\} .
\]
If we compare with the Riemann symbol of the standard hypergeometric
operator for $_2F_1(a,b,c;x)$
\[
\left\{ \begin{array}{ccc}
0&1&\infty\\
\hline
0&0&a\\
1-c&c-a-b&b\\
\end{array}\right\}
\]
we read off $a=(1-d+k)/2, b=(1+d+k)/2, c=1$, hence the
solution of our original operator is the square of
\[ (1-x)^{k/2} \cdot _2F_1((1-d+k)/2,(1+d+k)/2,1;x) .\]

\centerline{\bf \em  Corollary} 

{\em If $k$ and $d$ are rational numbers, then
the only the primes appearing in denominator of $k$ and $d$
appear in the denominator of the coefficients of the holomorphic
solution of $P$. After an appropriate rescaling, the coefficients
become integral.}\\

\newpage
\section*{\Large \bf \em Appendix C}
In this Appendix we list all essentially distinct Calabi-Yau operators 
of order four and degree two that are known to us. We also present some
further information.
To present the monodromy, it is sometimes convenient to use $u_i(x)=y_i(x)/(2\pi i)^3$ 
which we call the {\em scaled Frobenius basis}. The monodromy transformation around $0$ 
in this basis is given by the matrix
\[
\left(
\begin{array}{cccc}
1&0&0&0\\
1&1&0&0\\
1/2&1&1&0\\
1/6&1/2&1&1
\end{array}
\right) .
\]
In the cases with four singular points, the operator belongs to the main
component and the Riemann symbol of the operator has the form
\[ \left\{
\begin{array}{cccc} 
0&c_1&c_2&\infty\\
\hline
0&0&0&\alpha\\
0&1&1&\beta\\
0&1&1&2-\alpha\\
0&2&2&2-\beta\\
 \end{array} \right\}\]
where $c_1$, $c_2$ are solutions to the equation
\[ \Delta(x):=1+a x+f x^2 =0 \]
We will always assume that the exponents $\alpha, \beta$ are between $0$ and $1$. For each operator we give the monodromy transformations around the other singular points the with respect to this basis. This information was computed by {\sc J. Hofmann}, \cite{Hofmann}. 
A very common transformation is the {\em symplectic reflection} 
$v \mapsto v-<v,w>w$ in a vector $w=(a,b,c,d)$
represented by the matrix
\[I-\frac{1}{a}\left(\begin{array}{cccc}
-ad&ac&-ab&a^2\\
-bd&bc&-b^2&ab\\
-cd&c^2&-bc&ac\\ 
-d^2&cd&-bd&ad\\
\end{array} \right)\]
Such a tranformation is found at {\em conifold points}, where the exponents are $0\;\; 1\;\; 1\;\; 2$. 
For the conifold nearest to the origin the vector is of the form
\[ (a,0,c,d)\]
and in that case the {\em characteristic numbers $h^3,c_2h, c_3=\chi$} 
of the operator are determined by the reflection vector via
\[a=h^3,\;\;c=\frac{c_2h}{24},\;\;d=c_3 \lambda \]
where
\[ \lambda=\frac{\zeta(3)}{(2\pi i)^3} .\]
These numbers are formally attached to the differential operator, but have
the interpretation as characteristic numbers of the mirror manifold, see \cite{vS}.

{\bf \em Modular form information.} For conifold points appearing at rational values of the parameters one can determine the coeffients of a weight four modular
form as described in the thesis of {\sc K. Samol}, \cite{Samol}, \cite{SvS}.
We include here the name of modular form from $S_4(\Gamma_0(N)$ as it appears 
in the list of {\sc C. Meyer}, \cite{Meyer}. For example,  $(6/1)$ denotes the 
first (and only) weight four newform of level $6$,  which is the $\eta$-product
\[ f=(\eta(q)\eta(q^2)\eta(q^3)\eta(q^6))^2=q-2q^2-3q^3+4q^4+6q^5+\ldots \in S_4(\Gamma_0(6))  \]
For critical points at real quadratic irrationalities one expects the 
appearance of Hilbert modular forms; at imaginary quadratic irrationalities 
one would expect Bianchi modular forms. Another phenomenon that may occur is
the of a so-called {\em $K$-point}, see \cite{vS}. Here there are two single
logarithms appearing (two size two Jordan blocks), and one expects the 
appearance of a modular form from $S_3(\Gamma_0(N))$.

Below we list the $90$-degree rotated {\em extended Riemann-symbols}; 
each row starts with the singular value, followed by the four exponents 
at that point, then monodromy information and finally modular form information.
At the MUM point the monodromy is always standard, and we give the first
few instanton numbers as a substitute for modular form information.
If Hilbert-modular, Bianchi or weight three modular forms are expected, 
we indicate this with $h$, $b$ or $k$. In case the monodromy is of
finite order we put $-$.

The operators with three singular points which are on the main 
components can all be written in the form

\[\theta^4 + f x (-2\theta^4-4\theta^3+(\c(1-\c)+\a^2+\b^2-4)\theta^2+(\c(1-\c)+\a^2+\b^2-2)\theta+e+\]
\[+f^2(\theta + 1 - \a)(\theta + 1 - \b)(\theta + 1 + \b)(\theta + 1 + \a)\]
which has Riemann symbol
\[
\left\{
\begin{array}{ccc}
0&1/f&\infty\\
\hline
0&0&1-\a\\
0&\c&1-\b\\
0&1-\c&1+\b\\
0&1&1+\a\\
\end{array}
\right\}
\]
There is a single accessory parameter $e$.\\

\newpage

{\bf 2.1: $A\ast a =\#45$} \hfill $\chi=-120,\;\; c_2H = 72,\;\; H^3 =24,\;\; \dim |H| = 10$.

\[\theta^4-4 (2\theta+1)^2(7\theta^2+7\theta+2)-128(2\theta+1)^2(2\theta+3)^2\]

\[
\left\{\begin{array}{c|cccc||rrrr||c}
0&0&0&0&0&&&&&12,\;\;163,\;\; 3204\\
\frac{1}{128}&0&1&1&2&(24,&   0,&  3,&    -120 \lambda)&(64/5)\\[1mm]
-\frac{1}{16}&0&1&1&2&(48,& -24,& 10,&-3 - 240 \lambda)&(8/1)\\[1mm]
\infty&\frac{1}{2}&\frac{1}{2}&\frac{3}{2}&\frac{3}{2}&&&&&k\\
\end{array}
\right\}
\]
\vskip 10pt
{\bf 2.2: $B \ast a=\#15$} \hfill $\chi=-162,\;\; c_2H = 72,\;\; H^3 =18,\;\;\; \dim |H| = 9$.
\[\theta^4-3 x (3\theta+1)(3\theta+2)(7\theta^2+7\theta+2)-72 x^2 
(3\theta+1)(3\theta+2)(3\theta+4)(3\theta+5)\]
\[
\left\{\begin{array}{c|cccc||rrrr||c}
0&0&0&0&0&&&&&21,\;\; 480,\;\; 15894\\
\frac{1}{216}&0&1&1&2&(18, &  0, & 3,  &   -162 \lambda)&(54/2)\\[1mm]
-\frac{1}{27}&0&1&1&2&(36,& -18,& 9,& -3 - 324 \lambda)&(27/2)\\[1mm]
\infty&\frac{1}{3}&\frac{2}{3}&\frac{4}{3}&\frac{5}{3}&&&&&-\\
\end{array}
\right\}
\]
\vskip 10pt
{\bf 2.3: $C \ast a=\#68$} \hfill $\chi=-228,\;\; c_2H = 72,\;\; H^3 =12,\;\;\; \dim |H| = 8$.
\[\theta^4-4 x (4\theta+1)(4\theta+3)(7\theta^2+7\theta+2)-128 x (4\theta+1)(4\theta+3)(4\theta+5)(4\theta+7)\]
\[
\left\{\begin{array}{c|cccc||rrrr||c}
0&0&0&0&0&&&&&52,\;\; 2814,\;\; 220220\\
\frac{1}{512}&0&1&1&2&(12, &  0,& 3,   &  -228 \lambda)&(256/3)\\[1mm]
-\frac{1}{64}&0&1&1&2&(24,& -12,& 8,& -3 - 456 \lambda)&(32/3)\\[1mm]
\infty&\frac{1}{4}&\frac{3}{4}&\frac{5}{4}&\frac{7}{4}&&&&&-\\
\end{array}
\right\}
\]
\vskip 10pt
{\bf 2.4: $D \ast a=\#62$} \hfill $\chi=-336,\;\; c_2H = 72,\;\; H^3 =6,\;\;\; \dim |H| = 7$.

\[\theta^4-12 x (6\theta+1)(6\theta+5)(7\theta^2+7\theta+2)-1152x^2(6\theta+1)(6\theta+5)(6\theta+7)(6\theta+11)\]

\[
\left\{\begin{array}{c|cccc||rrrr||c}
0&0&0&0&0&&&&&372,\;\;136182,\;\;71562236\\
\frac{1}{3456}&0&1&1&2&(6,&   0,& 3,&     -366 \lambda)&(1728/16)\\[1mm]
-\frac{1}{432}&0&1&1&2&(12,& -6,& 7,& -3 - 732 \lambda)&(216/4)\\[1mm]
\infty&\frac{1}{6}&\frac{5}{6}&\frac{7}{6}&\frac{11}{6}&&&&&\\
\end{array}
\right\}
\]

{\bf 2.5: $A \ast b=\#25$}\hfill $\chi=-120,\;\; c_2H = 68,\;\; H^3 =20,\;\;\; \dim |H| = 9$.
\[\theta^4-4 x(2\theta+1)^2(11\theta^2+11\theta+3)-16x^2(2\theta+1)^2(2\theta+3)^2\]
\[
\left\{\begin{array}{c|cccc||rrrr||c}
0&0&0&0&0&&&&&20,\;\;277,\;\;8220\\
0.005636&0&1&1&2&(20,&   0,& 17/6,&     -120 \lambda)&h\\[1mm]
-0.693136&0&1&1&2&(80,&-40,&46/3&-13/3-480\lambda)&h\\[1mm]
\infty&\frac{1}{2}&\frac{1}{2}&\frac{3}{2}&\frac{3}{2}&&&&&k\\
\end{array}
\right\}
\]

\vskip 10pt

{\bf 2.6: $B \ast b=\#24$} \hfill $\chi=-150,\;\; c_2H = 66,\;\; H^3 =15,\;\;\; \dim |H| = 8$.
\[\theta^4-3 x (3\theta+1)(3\theta+2)(11\theta^2+11\theta+3)-9x^2(3\theta+1)(3\theta+2)(3\theta+4)(3\theta+5)\]

\[
\left\{\begin{array}{c|cccc||rrrr||c}
0&0&0&0&0&&&&&36,\;\;837,\;\;41421\\
0.003340&0&1&1&2&(15, &  0,& 11/4,  &   -150 \lambda)&h\\[1mm]
-0.410748&0&1&1&2&(12, &-6,& 7, &-3 - 732 \lambda)&h\\[1mm]
\infty&\frac{1}{3}&\frac{2}{3}&\frac{4}{3}&\frac{5}{3}&&&&&-\\
\end{array}
\right\}
\]
\vskip 10pt
{\bf 2.7: $C \ast b=\#51$}\hfill $\chi=-200,\;\; c_2H = 64,\;\; H^3 =10,\;\;\; \dim |H| = 7$.
\[\theta^4-4 x (4\theta+1)(4\theta+3)(11\theta^2+11\theta+3)-16x^2(4\theta+1)(4\theta+3)(4\theta+5)(4\theta+7)\]

\[
\left\{\begin{array}{c|cccc||rrrr||c}
0&0&0&0&0&&&&&92,\;\;5052,\;\;585396\\
0.0014090&0&1&1&2&(10, &  0,& 8/3, &    -200 \lambda)&h\\[1mm]
-0.173284&0&1&1&2&(40,&-20,&38/3,&-14/3-400\lambda)&h\\[1mm]
\infty&\frac{1}{4}&\frac{1}{2}&\frac{3}{2}&\frac{3}{2}&&&&&-\\
\end{array}
\right\}
\]
\vskip 10pt
{\bf 2.8 $D\ast b=\#63$}\hfill $\chi=-310,\;\; c_2H = 62,\;\; H^3 =5,\;\;\; \dim |H| = 6$.

\[\theta^4-12 x (6\theta+1)(6\theta+5)(11\theta^2+11\theta+3)-144x^2(6\theta+1)(6\theta+5)(6\theta+7)(6\theta+11)\]

\[
\left\{\begin{array}{c|cccc||rrrr||c}
0&0&0&0&0&&&&&684,\;\;253314,\;\;195638820\\
0.000208&0&1&1&2& (5, &  0,& 31/12,&  -310 \lambda)&h\\[1mm]
-0.025671&0&1&1&2&(20&-10&34/3 &-29/6-1240 \lambda )&h\\[1mm]
\infty&\frac{1}{6}&\frac{5}{6}&\frac{7}{6}&\frac{11}{6}&&&&&-\\
\end{array}
\right\}
\]

{\bf 2.9: $A \ast c=\#58$} \hfill $\chi=-112,\;\; c_2H = 72,\;\; H^3 =24,\;\;\; \dim |H| = 10$.
\[\theta^4-4 x(2\theta+1)^2(10\theta^2+10\theta+3)+144 x^2(2\theta+1)^2(2\theta+3)^2\]

\[
\left\{\begin{array}{c|cccc||rrrr||c}
0&0&0&0&0&&&&&16,\;\;142,\;\;11056/3;\\\
\frac{1}{144}&0&1&1&2&(24, &  0,& 3,   &  -112 \lambda)&(48/1)\\[1mm]
\frac{1}{16}&0&1&1&2&(72, &-24,& 9,& -3 - 336 \lambda)&(16/1)\\[1mm]
\infty&\frac{1}{2}&\frac{1}{2}&\frac{3}{2}&\frac{3}{2}&&&&&k\\
\end{array}
\right\}
\]
\vskip 10pt
{\bf 2.10: $B \ast c=\#70$}\hfill $\chi=-156,\;\; c_2H = 72,\;\; H^3 =18,\;\;\; \dim |H| = 9$.
\[\theta^4-3 x (3\theta+1)(3\theta+2)(10\theta^2+10\theta+3)-81x^2(3\theta+1)(3\theta+2)(3\theta+4)(3\theta+5)\]

\[
\left\{\begin{array}{c|cccc||rrrr||c}
0&0&0&0&0&&&&&27,\;\;n_2=432,\;\;n_3=18089\\
\frac{1}{243}&0&1&1&2&(18,  & 0, &3, &    -156 \lambda)&(243/1)\\[1mm]
\frac{1}{27}&0&1&1&2&(54,& -18,& 9,& -3 - 468 \lambda)&(27/1)\\[1mm]
\infty&\frac{1}{4}&\frac{3}{4}&\frac{5}{4}&\frac{7}{4}&&&&&-\\
\end{array}
\right\}
\]
\vskip 10pt
{\bf 2.11: $C \ast c=\#69$}\hfill $\chi=-224,\;\; c_2H = 72,\;\; H^3 =12,\;\;\; \dim |H| = 8$.
\[\theta^4-4 x (4\theta+1)(4\theta+3)(10\theta^2+10\theta+3)+144x^2(4\theta+1)(4\theta+3)(4\theta+5)(4\theta+7)\]

\[
\left\{\begin{array}{c|cccc||rrrr||c}
0&0&0&0&0&&&&&64,\;\;2616,\;\;246848;\;\\
\frac{1}{576}&0&1&1&2&(12, &  0,& 3, &    -224 \lambda)&(576/3)\\[1mm]
\frac{1}{64}&0&1&1&2&(36, &-12, &9, &-3 - 672 \lambda)&(64/3)\\[1mm]
\infty&\frac{1}{4}&\frac{3}{4}&\frac{5}{4}&\frac{7}{4}&&&&&-\\
\end{array}
\right\}
\]

\vskip 10pt
{\bf 2.12: $D \ast c=\#64$}\hfill $\chi=-364,\;\; c_2H = 72,\;\; H^3 =6,\;\; \dim |H| = 7$.
\[\theta^4-12 x (6\theta+1)(6\theta+5)(10\theta^2+10\theta+3)+1296 x^2(6\theta+1)(6\theta+5)(6\theta+7)(6\theta+11)\]
\[
\left\{\begin{array}{c|cccc||rrrr||c}
0&0&0&0&0&&&&&432,\;\;130842,\;\;78259376;\\
\frac{1}{3888}&0&1&1&2&(6, &  0,& 3,&     -364 \lambda)&(1944/5)\\[1mm]
\frac{1}{432}&0&1&1&2&(18,& -6,& 9,& -3 - 1092 \lambda)&(432/9)\\[1mm]
\infty&\frac{1}{4}&\frac{3}{4}&\frac{5}{4}&\frac{7}{4}&&&&&-\\
\end{array}
\right\}
\]

{\bf 2.13: $A \ast d=\#36$} \hfill $\chi=-88,\;\; c_2H = 80,\;\; H^3 =32,\;\;\; \dim |H| = 12$.
\[\theta^4-16 x(2\theta+1)^2(3\theta^2+3\theta+1)+512 x^2(2\theta+1)^2(2\theta+3)^2\]
\[
\left\{\begin{array}{c|cccc||rrrr||c}
0&0&0&0&0&&&&&16,\;\;42,\;\;1232\;\;\\
\frac{1}{128}&0&1&1&2&(32,  & 0, &10/3,&     -88 \lambda)&(64/4)\\[1mm]
\frac{1}{64}&0&1&1&2&(64,& -16,& 20/3,& -7/3 - 176 \lambda)&(32/2)\\[1mm]
\infty&\frac{1}{2}&\frac{1}{2}&\frac{3}{2}&\frac{3}{2}&&&&&k\\
\end{array}
\right\}
\]
\vskip 10pt
{\bf 2.14: $B \ast d=\#48$} \hfill $\chi=-162,\;\; c_2H = 84,\;\; H^3 =24,\;\;\; \dim |H| = 11$.
\[\theta^4-12 x (3\theta+1)(3\theta+2)(3\theta^2+3\theta+1)+288 x^2(3\theta+1)(3\theta+2)(3\theta+4)(3\theta+5)\]

\[
\left\{\begin{array}{c|cccc||rrrr||c}
0&0&0&0&0&&&&&4,\;\;291/2,\;\;5832\\
\frac{1}{216}&0&1&1&2&(24, &  0,& 7/2,&     -162 \lambda)&(9/1)\\[1mm]
\frac{1}{108}&0&1&1&2&(48, &-12,& 7,& -9/4 - 324 \lambda)&(108/4)\\[1mm]
\infty&\frac{1}{3}&\frac{2}{3}&\frac{4}{3}&\frac{5}{3}&&&&&-\\
\end{array}
\right\}
\]

\vskip 10pt
{\bf 2.15: $C \ast d=\#38$}\hfill $\chi=-268,\;\; c_2H = 88,\;\; H^3 =16,\;\;\; \dim |H| = 10$.
\[\theta^4-16 x (4\theta+1)(4\theta+3)(3\theta^2+3\theta+1)+512x^2(4\theta+1)(4\theta+3)(4\theta+5)(4\theta+7)\]

\[
\left\{\begin{array}{c|cccc||rrrr||c}
0&0&0&0&0&&&&&48,\;\;998,\;\;73328\;\\
\frac{1}{512}&0&1&1&2&(16, &  0,& 11/3, &    -268 \lambda)&(256/1)\\[1mm]
\frac{1}{256}&0&1&1&2&(32,& -8,& 22/3,& -13/6 - 536 \lambda)&(128/4)\\[1mm]
\infty&\frac{1}{4}&\frac{3}{4}&\frac{5}{4}&\frac{7}{4}&&&&&-\\
\end{array}
\right\}
\]

\vskip 10pt
{\bf 2.16: $D \ast d=\#65$}\hfill $\chi=-470,\;\; c_2H =92,\;\; H^3 =8,\;\; \dim |H| = 9$.
\[\theta^4-48 x (6\theta+1)(6\theta+5)(3\theta^2+3\theta+1)+4608 x^2(6\theta+1)(6\theta+5)(6\theta+7)(6\theta+11)\]

\[
\left\{\begin{array}{c|cccc||rrrr||c}
0&0&0&0&0&&&&&240,\;\;57102,\;\;19105840;\\
\frac{1}{3456}&0&1&1&2&(8,&   0,& 23/6, &    -470 \lambda)&(576/8)\\[1mm]
\frac{1}{1728}&0&1&1&2&(16,& -4,& 23/3, &-25/12 - 940 \lambda)&(864/3)\\[1mm]
\infty&\frac{1}{6}&\frac{5}{6}&\frac{7}{6}&\frac{11}{6}&&&&&-\\
\end{array}
\right\}
\]

{\bf 2.17: $A \ast e=\#111$}
\[\theta^4-16x(2\theta+1)^2(8\theta^2+8\theta+3)+2^{12} x^2(2\theta+1)^2(2\theta+3)^2\]
\[
\left\{
\begin{array}{c|cccc||c}
0&0&0&0&0&32, -96, 1440 \\
1/256&0&1/2&1/2&1&m\\
\infty&1/2&1/2&3/2&3/2&k\\
\end{array}
\right\}
\]
\vskip 20pt
{\bf 2.18: $B \ast e =\#110$}

\[\theta^4-12x(3\theta+1)(3\theta+2)(8\theta^2+8\theta+3)+2^83^2x^2(3\theta+1)(3\theta+2)(3\theta+4)(3\theta+5)\]
\[
\left\{
\begin{array}{c|cccc||c}
0&0&0&0&0&36,-144,8076\\
1/432&0&1/2&1/2&1&m\\
\infty&1/3&2/3&4/3&5/3&-\\
\end{array}
\right\}
\]
\vskip 20pt
{\bf 2.xx19: $C \ast e \sim 1.3=\#3$}
\[\theta^4-16x(4\theta+1)(4\theta+3)(8\theta^2+8\theta+3)+2^{12}x^2(4\theta+1)(4\theta+3)(4\theta+5)
(4\theta+7)\]
\[
\left\{
\begin{array}{c|cccc||c}
0&0&0&0&0&32,608,26016\\
1/1024&0&1/2&1/2&1&m\\
\infty&1/4&3/4&5/4&7/4&-\\
\end{array}
\right\}
\]
\vskip 20pt

{\bf 2.19: $D \ast e=\#112$}
\[\theta^4-48x(6\theta+1)(6\theta+5)(8\theta^2+8\theta+3)+2^{12} 3^2 x^2(6\theta+1)(6\theta+5)(6\theta+7)(6\theta+11)\]

\[
\left\{
\begin{array}{c|cccc||c}
0&0&0&0&0&-288, 162504, -96055968\\
1/6912&0&1/2&1/2&1& m\\
\infty&1/6&5/6&7/6&11/6&-\\
\end{array}
\right\}
\]

{\bf 2.20: $A \ast f=\#133$}
\[\theta^4-12 x(2\theta+1)^2(3\theta^2+3\theta+1)+432x^2(2\theta+1)^2(2\theta+3)^2\]
\[
\left\{\begin{array}{c|cccc||rrrr||c}
0&0&0&0&0&&&&&12,\;\;-42,\;\;-3284/3\\
\alpha&0&1&1&2&(36, &  -6,& 4, & -5/6   -120 \lambda)&b\\[1mm]
\beta&0&1&1&2&(36, &  6,& 4,& -5/6 - 120 \lambda)&b\\[1mm]
\infty&\frac{1}{2}&\frac{1}{2}&\frac{3}{2}&\frac{3}{2}&&&&&k\\
\end{array}
\right\}
\]
\vskip 10pt
{\bf 2.21: $B \ast f=\#134$}
\[\theta^4-9 x (3\theta+1)(3\theta+2)(3\theta^2+3\theta+1)+243 x^2(3\theta+1)(3\theta+2)(3\theta+4)(3\theta+5)\]

\[
\left\{\begin{array}{c|cccc||rrrr||c}
0&0&0&0&0&&&&&18,\;\;-207/2,\;\;-52177\\
\alpha&0&1&1&2&(27,& -9/2,& 33/8,& -13/16    -198 \lambda)&b\\[1mm]
\beta&0&1&1&2& (27,&  9/2,& 33/8,& -13/16 - 198 \lambda)&b\\[1mm]
\infty&\frac{1}{3}&\frac{2}{3}&\frac{4}{3}&\frac{5}{3}&&&&&-\\
\end{array}
\right\}
\]
\vskip 10pt

{\bf 2.22: $C \ast f=\#135$}
\[\theta^4-12 x (4\theta+1)(4\theta+3)(3\theta^2+3\theta+1)+432x^2(4\theta+1)(4\theta+3)(4\theta+5)(4\theta+7)\]
\[
\left\{\begin{array}{c|cccc||rrrr||c}
0&0&0&0&0&&&&&36,\;\;-477,\;\;-206716/3\\
\alpha&0&1&1&2&(18,  & -3, &17/4,&  -19/24   -312 \lambda)&b\\[1mm]
\beta&0&1&1&2&(18, &3,& 17/4,& 19/24 - 312 \lambda)&b\\[1mm]
\infty&\frac{1}{4}&\frac{3}{4}&\frac{5}{4}&\frac{7}{4}&&&&&-\\
\end{array}
\right\}
\]
\vskip 10pt
{\bf 2.23: $D \ast f =\#136$}
\[\theta^4-36 x (6\theta+1)(6\theta+5)(3\theta^2+3\theta+1)+3888 x^2(6\theta+1)(6\theta+5)(6\theta+7)(6\theta+11)\]
\[
\left\{\begin{array}{c|cccc||rrrr||c}
0&0&0&0&0&&&&&180,\;\;-15615,\;\;-21847076\;\\
\alpha&0&1&1&2&(9, &  -3/2, &35/8,&  -37/48 -534 \lambda)&b\\[1mm]
\beta&0&1&1&2&(9,& 3/2, &35/8,& -37/48 - 534 \lambda)&b\\[1mm]
\infty&\frac{1}{6}&\frac{5}{6}&\frac{7}{6}&\frac{11}{6}&&&&&-\\
\end{array}
\right\}
\]

{\bf 2.24: $A \ast g=\#137$}\hfill $\chi=-16,\;\; c_2H =96,\;\; H^3 =48,\;\;\; \dim |H| = 16$.
\[\theta^4-4 x(2\theta+1)^2(17\theta^2+17\theta+6)+1152 x^2(2\theta+1)^2(2\theta+3)^2\]

\[
\left\{\begin{array}{c|cccc||rrrr||c}
0&0&0&0&0&&&&&20,\;\;2,\;\;1684/3\\
\frac{1}{144}&0&1&1&2&(48,  & 0,& 4,&     -16 \lambda)&(24/1)\\[1mm]
\frac{1}{128}&0&1&1&2&(72,& -12,& 6,& -2 - 24 \lambda)&(64/1)\\[1mm]
\infty&\frac{1}{2}&\frac{1}{2}&\frac{3}{2}&\frac{3}{2}&&&&&k\\
\end{array}
\right\}
\]
\vskip 10pt
{\bf 2.25: $B \ast g=\#138$}\hfill $\chi=-156,\;\; c_2H =108,\;\; H^3 =36,\;\;\; \dim |H| = 15$.
\[\theta^4-3 x (3\theta+1)(3\theta+2)(17\theta^2+17\theta+6)+648 x^2(3\theta+1)(3\theta+2)(3\theta+4)(3\theta+5)\]

\[
\left\{\begin{array}{c|cccc||rrrr||c}
0&0&0&0&0&&&&&27,\;\;189/4,\;\;2618\\
\frac{1}{243}&0&1&1&2&(36,  & 0,& 9/2,&     -156 \lambda)&(243/2)\\[1mm]
\frac{1}{216}&0&1&1&2&(54,& -9, &27/4, &-15/8 - 234 \lambda)&(54/4)\\[1mm]
\infty&\frac{1}{3}&\frac{2}{3}&\frac{4}{3}&\frac{5}{3}&&&&&-\\
\end{array}
\right\}
\]
\vskip 10pt
{\bf 2.26: $C \ast g=\#139$}\hfill $\chi=-344,\;\; c_2H =120,\;\; H^3 =24,\;\;\; \dim |H| = 14$.
\[\theta^4-4 x (4\theta+1)(4\theta+3)(17\theta^2+17\theta+6)+1152 x^2(4\theta+1)(4\theta+3)(4\theta+5)(4\theta+7)\]

\[
\left\{\begin{array}{c|cccc||rrrr||c}
0&0&0&0&0&&&&&44,\;\;607,\;\;22500\;\\
\frac{1}{576}&0&1&1&2&(24,  & 0,& 5, &    -344 \lambda)&(288/10)\\[1mm]
\frac{1}{512}&0&1&1&2&(36, &-6,& 15/2,& -7/2 - 516 \lambda)&(256/4)\\[1mm]
\infty&\frac{1}{4}&\frac{3}{4}&\frac{5}{4}&\frac{7}{4}&&&&&-\\
\end{array}
\right\}
\]
\vskip 10pt

{\bf 2.27: $D \ast g=\#140$}\hfill $\chi=-676,\;\; c_2H =132,\;\; H^3 =12,\;\;\dim |H| = 13$.
\[\theta^4-12 x (6\theta+1)(6\theta+5)(17\theta^2+17\theta+6)+10368 x^2(6\theta+1)(6\theta+5)(6\theta+7)(6\theta+11)\]

\[
\left\{\begin{array}{c|cccc||rrrr||c}
0&0&0&0&0&&&&&108,\;\;54135,\;\;-494556\\
\frac{1}{3888}&0&1&1&2&(12,&   0, &11/2,&     -676 \lambda)&(1944/6)\\[1mm]
\frac{1}{3456}&0&1&1&2&(18,& -3,& 33/4,& -13/8 -1014 \lambda)&(1728/15)\\[1mm]
\infty&\frac{1}{6}&\frac{5}{6}&\frac{7}{6}&\frac{11}{6}&&&&&-\\
\end{array}
\right\}
\]

{\bf 2.28: $A \ast h=\#141$}
\[\theta^4-12x(2\theta+1)^2(18\theta^2+18\theta+7)+2^43^6x^2(2\theta+1)^2(2\theta+3)^2\]
\[
\left\{
\begin{array}{c|cccc||c}
0&0&0&0&0&48, -438, 2864\\
1/432&0&1/3&2/3&1&-\\
\infty&1/2&1/2&3/2&3/2&k\\
\end{array}
\right\}
\]
\vskip 20pt

{\bf 2.29: $B \ast h=\#142$} 
\[\theta^4-9x(3\theta+1)(3\theta+2)(18\theta^2+18\theta+7)+3^8 x^2(3\theta+1)(3\theta+2)(3\theta+4)(3\theta+5)\]
\[
\left\{
\begin{array}{c|cccc||c}
0&0&0&0&0&45, -3465/4, 27735\\
1/729&0&1/3&2/3&1&-\\
\infty&1/3&2/3&4/3&5/3&-\\
\end{array}
\right\}
\]
\vskip 20pt
{\bf 2.xx30: $C \ast h$} 
\[\theta^4-12x(4\theta+1)(4\theta+3)(18\theta^2+18\theta+7)+2^{4}3^6x^2(4\theta+1)(4\theta+3)(4\theta+5)(4\theta+7)\]
\[
\left\{
\begin{array}{c|cccc||c}
0&0&0&0&0&0,0,0\\
1/1728&0&1/3&2/3&1&-\\
\infty&1/4&3/4&5/4&7/4&-\\
\end{array}
\right\}
\]
\vskip 20pt
{\bf 2.30: $D \ast h=\#143$}
\[\theta^4-36x(6\theta+1)(6\theta+5)(18\theta^2+18\theta+7)+2^{4} 3^8 x^2(6\theta+1)^(6\theta+5)(6\theta+7)(6\theta+11)\]
\[
\left\{
\begin{array}{c|cccc||c}
0&0&0&0&0& -1008, 499086, -607849200\\
1/11664&0&1/3&2/3&1&-\\
\infty&1/6&5/6&7/6&11/6&-\\
\end{array}
\right\}
\]

{\bf 2.31: $A \ast i$}
\[\theta^4-16x(2\theta+1)^2(32\theta^2+32\theta+13)+2^{16}x^2(2\theta+1)^2(2\theta+3)^2\]
\[
\left\{
\begin{array}{c|cccc||c}
0&0&0&0&0&96, -3560, -12064\\
1/1024&0&1/4&3/4&1&-\\
\infty&1/2&1/2&3/2&3/2&k\\
\end{array}
\right\}
\]

\vskip 20pt
{\bf 2.32: $B \ast i$} 
\[\theta^4-12x(3\theta+1)(3\theta+2)(32\theta^2+32\theta+13)+2^{12}3^2 x^2(3\theta+1)(3\theta+2)(3\theta+4)(3\theta+5)\]
\[
\left\{
\begin{array}{c|cccc||c}
0&0&0&0&0& 60, -7635, 307860 \\
1/1728&0&1/4&3/4&1&-\\
\infty&1/3&2/3&4/3&5/3&-\\
\end{array}
\right\}
\]
\vskip 20pt
{\bf 2.33: $C \ast i$} 
\[\theta^4-16x(4\theta+1)(4\theta+3)(32\theta^2+32\theta+13)+2^{4}3^6x^2(4\theta+1)(4\theta+3)(4\theta+5)(4\theta+7)\]
\[
\left\{
\begin{array}{c|cccc||c}
0&0&0&0&0&-160, -6920, -539680\\
1/4096&0&1/4&3/4&1&-\\
\infty&1/4&3/4&5/4&7/4&-\\
\end{array}
\right\}
\]
\vskip 20pt
{\bf 2.34: $D \ast i$}
\[\theta^4-48x(6\theta+1)(6\theta+5)(32\theta^2+32\theta+13)+2^{16} 3^2 x^2(6\theta+1)^(6\theta+5)(6\theta+7)(6\theta+11)\]
\[
\left\{
\begin{array}{c|cccc||c}
0&0&0&0&0&-3936, 3550992, -10892932064\\
1/27648&0&1/4&3/4&1&-\\
\infty&1/6&5/6&7/6&11/6&-\\
\end{array}
\right\}
\]
\vskip20pt
{\bf 2.35: $A \ast j$}
\[\theta^4-48x(2\theta+1)^2(72\theta^2+72\theta+31)+2^{12}3^6x^2(2\theta+1)^2(2\theta+3)^2\]
\[
\left\{
\begin{array}{c|cccc||c}
0&0&0&0&0&480, -226968, -16034720\\
1/6912&0&1/6&5/6&1&-\\
\infty&1/2&1/2&3/2&3/2&k\\
\end{array}
\right\}
\]
\vskip 20pt
{\bf 2.36: $B \ast j$}
\[\theta^4-36x(3\theta+1)(3\theta+2)(72\theta^2+72\theta+31)+2^{8}3^8 x^2(3\theta+1)(3\theta+2)(3\theta+4)(3\theta+5)\]
\[
\left\{
\begin{array}{c|cccc||c}
0&0&0&0&0&-36, -486279, 128217204\\
1/11664&0&1/6&5/6&1&-\\
\infty&1/3&2/3&4/3&5/3&-\\
\end{array}
\right\}
\]
\vskip 20pt
{\bf 2.37: $C \ast j$} 
\[\theta^4-48x(4\theta+1)(4\theta+3)(72\theta^2+72\theta+31)+2^{12}3^6x^2(4\theta+1)(4\theta+3)(4\theta+5)(4\theta+7)\]
\[
\left\{
\begin{array}{c|cccc||c}
0&0&0&0&0&-2592, -307800, 81451104\\
1/27648&0&1/6&5/6&1&-\\
\infty&1/4&3/4&5/4&7/4&-\\
\end{array}
\right\}
\]
\vskip 20pt
{\bf 2.38: $D \ast j$}
\[\theta^4-144x(6\theta+1)(6\theta+5)(72\theta^2+72\theta+31)+2^{12} 3^8 x^2(6\theta+1)^(6\theta+5)(6\theta+7)(6\theta+11)\]
\[
\left\{
\begin{array}{c|cccc||c}
0&0&0&0&0&-41184, 251271360, -5124430612320\\
1/186624&0&1/6&5/6&1&-\\
\infty&1/6&5/6&7/6&11/6&-\\
\end{array}
\right\}
\]
\vskip 20pt
{\bf 2.52: $I \ast \alpha =\#16$}\hfill $\chi=-128,\;\; c_2H =96,\;\; H^3 =48,\;\;\dim |H| = 16$.
\[\theta^4-4x (2\theta+1)^2(5\theta^2+5\theta+2)+256 x^2 (2\theta)(\theta+1)^2(2\theta+3)\]
\[
\left\{\begin{array}{c|cccc||rrrr||c}
0&0&0&0&0&&&&&4,\;\;20,\;\;644/3\\
\frac{1}{64}&0&1&1&2&(48,&   0,& 4,&     -128 \lambda)&(6/1)\\[1mm]
\frac{1}{16}&0&1&1&2&(192,& -48, &16,& -4 -512 \lambda)&(12/1)\\[1mm]
\infty&\frac{1}{2}&1&1&\frac{3}{2}&&&&&m\\
\end{array}
\right\}
\]
\vskip20pt
{\bf 2.53: $I \ast \gamma = \#29$}\hfill $\chi=-116,\;\; c_2H =72,\;\; H^3 =24,\;\;\dim |H| = 10$.
\[\theta^4-2x (2\theta+1)^2(17\theta^2+17\theta+5)+4 x^2 (2\theta+1)(\theta+1)^2(2\theta+3)\]
\[
\left\{\begin{array}{c|cccc||rrrr||c}
0&0&0&0&0&&&&&14,\;\;303/2,\;\;10424/3\;\\
0.00737&0&1&1&2&(24, &  0,& 3,&     -116 \lambda)&-\\[1mm]
8.49263&0&1&1&2&(600,&-240,&75&-20-2900 \lambda)&-\\[1mm]
\infty&\frac{1}{2}&1&1&\frac{3}{2}& & & & &m\\
\end{array}
\right\}
\]
\vskip 20pt
{\bf 2.54: $I \ast \delta =\#41$}\hfill $\chi=-116,\;\; c_2H =72,\;\; H^3 =24,\;\; \dim |H| = 10$.
\[\theta^4-2x (2\theta+1)^2(7\theta^2+7\theta+3)+324x^2 (2\theta+1)(\theta+1)^2(2\theta+3)\]

\[
\left\{\begin{array}{c|cccc||rrrr||c}
0&0&0&0&0&&&&&14,\;\;303/2,\;\;10424/3\\
*&0&1&1&2&(72,  & -12,& 6,&  -1   -180 \lambda)&b\\[1mm]
*&0&1&1&2&(72, &12,& 6,& 1 -180 \lambda)&b\\[1mm]
\infty&\frac{1}{2}&1&1&\frac{3}{2}&&&&&m\\
\end{array}
\right\}
\]
\vskip 20pt
{\bf 2.55: $I \ast \epsilon =\#42$}\hfill $\chi=-116,\;\; c_2H =80,\;\; H^3 =32,\;\; \dim |H| = 12. $
\[\theta^4-8 x (2\theta+1)^2(3\theta^2+3\theta+1)+64 x^2 (2\theta+1)(\theta+1)^2(2\theta+3)\]
\[
\left\{\begin{array}{c|cccc||rrrr||c}
0&0&0&0&0&&&&&8,\;\;\;63,\;\;\;1000\\
\alpha&0&1&1&2&(32, &  0,& 10/3,& -116 \lambda)&-\\[1mm]
\beta&0&1&1&2&(288, &  -96,& 30,& -8 -1044 \lambda)&-\\[1mm]
\infty&\frac{1}{2}&1&1&\frac{3}{2}&\\
\end{array}
\right\}
\]
\vskip 20pt
{\bf 2.56: $I \ast \zeta = \#185$} \hfill $\chi=-120,\;\; c_2H =84,\;\; H^3 =36,\;\;\; \dim |H| = 13 $.
\[\theta^4-6 x (2\theta+1)^2(3\theta^2+3\theta+1)-108 x^2 (2\theta+1)(\theta+1)^2(2\theta+3)\]

\[
\left\{\begin{array}{c|cccc||rrrr||c}
0&0&0&0&0&&&&&6,\;\;\;93/2,\;\;\;608\\
\alpha&0&1&1&2&(36,  & 0,& 7/2,& -120 \lambda)&h\\[1mm]
\beta&0&1&1&2&(144, &  -72,& 26,& -7 -480 \lambda)&h\\[1mm]
\infty&\frac{1}{2}&1&1&\frac{3}{2}&&&&&m\\
\end{array}
\right\}
\]
\vskip 20pt
{\bf 2.57: $I \ast \eta =\#184$} 
\[\theta^4-2 x (2\theta+1)^2(11\theta^2+11\theta+5)+500 x^2 (2\theta+1)(\theta+1)^2(2\theta+3)\]

\[
\left\{\begin{array}{c|cccc||rrrr||c}
0&0&0&0&0&&&&&2,\;\;\;4,\;\;\;-8\\\
\alpha&0&1&1&2&(100,&-10,& 20/3,&-5/6 -200 \lambda)&b\\[1mm]
\beta&0&1&1&2&(100, & 10, &20/3, &5/6 -200 \lambda)&b\\[1mm]
\infty&\frac{1}{2}&1&1&\frac{3}{2}&&&&&m\\
\end{array}
\right\}
\]
\vskip 20pt
{\bf 2.58: $I \ast \iota=4*$}
\[\theta^4 -6x(2\theta+1)^2(9\theta^2+9\theta+5)+ 2916x^2(2\theta+1)(\theta+1)^2(2\theta+3)\]
\[
\left\{
\begin{array}{c|cccc||c}
0&0&0&0&0&-6,-6,-104\\
1/108&0&1/6&5/6&1&-\\
\infty&1/2&1&1&3/2&m\\
\end{array}
\right\}
\]
\vskip 20pt
{\bf 2.xx59: $I \ast \theta \sim 2.17$.}

\[{\theta}^{4} -16\, x \left( 2\,\theta+1 \right) ^{2} \left( 8\,{\theta}^{2}+8\,\theta+5 \right) + 16384\, x^2 \left( 2\,\theta+1 \right)  \left( \theta+1
 \right) ^{2} \left( 2\,\theta+3 \right) \]
\[
\left\{
\begin{array}{c|cccc||c}
0&0&0&0&0&-32, -88, -1440\\
1/256&0&0&1&1&k\\
\infty&1/2&1&1&3/2&m\\
\end{array}
\right\}
\]
\vskip 20pt
{\bf 2.59: $I \ast \kappa$}
\[ \theta^4 -48 (2\theta+1)^2(18\theta^2+18\theta+13)+ 746496 x^2(2\theta+1)(\theta+1)^2(2\theta+3)
\]

\[
\left\{
\begin{array}{c|cccc||c}
0&0&0&0&0&-384, -1356, -164736\\
1/1728&-1/6&0&1&7/6&-\\
\infty&1/2&1&1&3/2&m\\
\end{array}
\right\}
\]
\vskip 20pt
{\bf 2.60: $=I \ast Sporadic 1=\#18$} \hfill $\chi=-128,\;\; c_2H =88,\;\; H^3 =40,\;\;\; \dim |H| = 14$.
\[\theta^4-4 x (2\theta+1)^2(3\theta^2+3\theta+1)-16 x^2 (2\theta+1)(4\theta+3)(4\theta+5)(2\theta+3)\]
\[
\left\{\begin{array}{c|cccc||rrrr||c}
0&0&0&0&0&&&&&4,\;\;39,\;\;364\\
\frac{1}{64}&0&1&1&2&(40,&0, &11/3,&-128 \lambda)&(80/3)\\[1mm]
-\frac{1}{16}&0&1&1&2&(80, & -40,& 46/3,& -13/3 -256 \lambda)&(40/2)\\[1mm]
\infty&\frac{1}{2}&\frac{3}{4}&\frac{5}{4}&\frac{3}{2}&&&&&-\\
\end{array}
\right\}
\]
\vskip 20pt
{\bf 2.61: $I \ast Sporadic 2=\#26$} \hfill  $\chi=-116,\;\; c_2H =76,\;\; H^3 =28,\;\;\; \dim |H| =11 $

\[\theta^4- 2x (2\theta+1)^2(13\theta^2+13\theta+4)-12 x^2 (2\theta+1)(3\theta+2)(3\theta+4)(2\theta+3)\]
\[
\left\{\begin{array}{c|cccc||rrrr||c}
0&0&0&0&0&&&&&10,\;\;191/2,\;\;\;1724\\
\frac{1}{108}&0&1&1&2&(28,&0, &19/6,&-116 \lambda)&(252/3)\\[1mm]
-\frac{1}{4}&0&1&1&2&(112,&  -56, &62/3,& -17/3 -464 \lambda)&(28/1)\\[1mm]
\infty&\frac{1}{2}&\frac{2}{3}&\frac{4}{3}&\frac{3}{2}&&&&&-\\
\end{array}
\right\}
\]
\vskip 20pt
{\bf 2.62: $=\#28$} \hfill $\chi=-96,\;\; c_2H =84,\;\; H^3 =42,\;\;\; \dim |H|= 14$.
\[\theta^4- x (65\theta^4+130\theta^3+105\theta^2+40\theta+6)+4 x^2 (4\theta+3)(\theta+1)^2(4\theta+5)\]

\[
\left\{\begin{array}{c|cccc||rrrr||c}
0&0&0&0&0&&&&&5,\;\;\;28,\;\;\;312\\
\frac{1}{64}&0&1&1&2&(42,&0,& 7/2,& -96 \lambda)&(14/2)\\[1mm]
1&0&1&1&2&(756,&  -252,& 63,& -15 -1728 \lambda)&(7/1)\\[1mm]
\infty&\frac{3}{4}&1&1&\frac{5}{4}&&&&&m\\
\end{array}
\right\}
\]

\[a_n=\sum_{j,k} {n \choose k}^2{n \choose j}^2{k+j \choose n}^2\]
The operator comes from the mirror symmetry of the
intersection of six general hyperplanes of the Grassmanian $G(3,6)$, 
Pl\"ucker-embedded in ${\bf P}^{19}$, \cite{BCKvS}.\\
\vskip 20pt

{\bf 2.63: $\#84$}
\[\theta^4-4x(32\theta^4+64\theta^3+63\theta^2+31\theta+6)+256 x^2(4\theta+3)(\theta+1)^2(4\theta+5)\]
 \[
\left\{
\begin{array}{c|cccc||c}
0&0&0&0&0& -4,\;\; -11,\;\;-44,\ldots\\
1/64&0&0&1&1&k\\
\infty&3/4&1&1&5/4&m\\
\end{array}
\right\}
\]
\[A_n=\sum_k {n \choose k}^2{2n \choose 2k}^{-1} \frac{(4k!)}{(2k!)k!^2}\frac{(4n-4k)!}{(2n-2k)!(n-k)!^2}\]
\vskip 20pt
{\bf 2.64: $\#182$} \hfill $\chi=-96,\;\; c_2H =132,\;\; H^3 =132,\;\;\; \dim |H| = $
\[\theta^4- x (43\theta^4+86\theta^3+77\theta^2+34\theta+6)+
12 x^2 (6\theta+5)(\theta+1)^2(6\theta+7)\]

\[ 
\left\{\begin{array}{c|cccc||rrrr||c}
0&0&0&0&0&&&&&1,\;\;7/4,\;\;7\\
\frac{1}{16}&0&1&1&2&(132,&0, &11/2,& -96 \lambda)&(22/3)\\[1mm]
\frac{1}{27}&0&1&1&2&(396,&-66,&33/2,&-13/4 -288 \lambda)&(33/2)\\[1mm]
\infty&\frac{5}{6}&1&1&\frac{7}{6}&&&&&m\\
\end{array}
\right\}
\]
This operator was found by a brute-force search.  We do not know an
explicit formula for the coefficient $A_n$.\\ 

\vskip 20pt
{\bf 2.65 =$I \ast Sporadic 3=\#183$} \hfill $\chi=-128,\;\;c_2H =120,\;\;H^3 =72,\;\;\dim |H| = 22$.
\[\theta^4-4x (2\theta+1)^2(7\theta^2+7\theta+3)+48 x^2(2\theta+1)(4\theta+3)(4\theta+5)(2\theta+3)\]

\[
\left\{\begin{array}{c|cccc||rrrr||c}
0&0&0&0&0&&&&&4,\;\;\;7,\;\;\;556/9\\
\frac{1}{64}&0&1&1&2&(72,&0,& 5,& -128 \lambda)&(16/1)\\[1mm]
\frac{1}{48}&0&1&1&2&(144,&  -24,& 10, &-7/3 -256 \lambda)&(72/1)\\[1mm]
\infty&\frac{1}{2}&\frac{3}{4}&\frac{5}{4}&\frac{3}{2}&&&&&-\\
\end{array}
\right\}
\]
\vskip 20pt
{\bf 2.66} Reducible operator (does not really count).

\[\theta^4 -12x(6\theta+1)(6\theta+5)(2\theta^2+2\theta+1),
+144x^2(6\theta+1)(6\theta+5)(6\theta+7)(6\theta +11)\]
\[
\left\{
\begin{array}{c|cccc||c}
0&0&0&0&0& -192,\;\;4182,\;\;-229568\\
1/432&0&0&1&1&k\\
\infty&1/6&5/6&7/6&11/6&-\\
\end{array}
\right\}
\]
This operator is the square of the second order hypergeometric 
operator $D$.\\

\vskip 20pt
{\bf 2.67: $\sim \#245$} (Bogner 1)

\[\theta^4-2x(108\theta^4+198\theta^3+183\theta^2+84\theta+15)+36x^2(3\theta+2)^2(6\theta+7)^2\]
\[
\left\{
\begin{array}{c|cccc||c}
0&0&0&0&0& -6,\;\;-33 ,\;\;-170\\
1/108&0&1/6&1&7/6&-\\
\infty&2/3&2/3&7/6&7/6&-\\
\end{array}
\right\}
\]
We do not know a formula for $A_n$.\\
\vskip 20pt
{\bf 2.68: $\sim \#406$} (Bogner 2)
\[\theta^4-4x(128 \theta^4+ 224 \theta^3+197 \theta^2+85 \theta + 14)+128x^2(2\theta+1)(4\theta+5)(8\theta+5)(8\theta+9)\]
\[
\left\{
\begin{array}{c|cccc||c}
0&0&0&0&0& -12\;\;,-186\;\;,-1668,\;\;\\
1/256&0&1/4&1&5/4&-\\
\infty&1/2&5/8&9/8&5/4&-\\
\end{array}
\right\}
\]
We do not know a formula for $A_n$.\\
\vskip 20pt
{\bf 2.69: $=\#205$} \hfill $\chi=-128,\;\; c_2H =160,\;\; H^3 =160,\;\;\; \dim |H| =40$.
\[\theta^4- x (59\theta+118\theta^3+105\theta^2+46\theta+8)+96 x^2 (3\theta+2)(\theta+1)^2(3\theta+4)\]
\[
\left\{\begin{array}{c|cccc||rrrr||c}
0&0&0&0&0&&&&&1,\;\;7/4,\;\;5\\
\frac{1}{32}&0&1&1&2&(160,&0,&20/3, & -128\lambda)&(5/1)\\[1mm]
\frac{1}{27}&0&1&1&2&(320,&-40,& 40/3,&-7/3 -256 \lambda)&(15/2)\\[1mm]
\infty&\frac{2}{3}&1&1&\frac{4}{3}&&&&&m\\
\end{array}
\right\}
\]

\[A_n =4 \sum_k^{\lfloor n/4 \rfloor} \frac{n-2k}{3n-4k}{n \choose k}{2k \choose k}{2n-2k \choose n-k}{3n-4k \choose 2n}\]
\vskip 20pt
{\bf 2.70: $\sim \#255$} (Bogner 3) 
\[\theta^4-4x(128\theta^4+160\theta^3+125\theta^2+45\theta+6)+128 x^2(8\theta + 7) (2 \theta + 1) (4 \theta + 3) (8 \theta + 3)\]
\[
\left\{
\begin{array}{c|cccc||c}
0&0&0&0&0& 20,\;\; 290,\;\;28820/3\\
1/256&0&3/4&1&7/4&-\\
\infty&3/8&1/2&3/4&7/8&-\\
\end{array}
\right\}
\]
We do not know a formula for $A_n$.\\

\vskip 20pt
{\bf \em The 14 tilde operators}
There are $14$ exponents $(\alpha_1,\alpha_2,\alpha_3,\alpha_4)$ with
\[\alpha_1 \le \alpha_2 \le \alpha_3 \le \alpha_4,\;\;\alpha_1+\alpha_4=\alpha_2+\alpha_4=1,\;\] 
for which the hypergeometric operator, scaled by $N$,
\[\theta^4-Nx(\theta+\alpha_1)(\theta+\alpha_2)(\theta+\alpha_3)(\theta+\alpha_4)\]
is a Calabi-Yau operator, \cite{Alm1}. Corresponding to these, there are also 
$14$ hypergeometric {\em fifth order} Calabi-Yau operator
\[\theta^5-4Nx(\theta+\alpha_1)(\theta+\alpha_2)(\theta+\frac{1}{2})(\theta+\alpha_3)(\theta+\alpha_4)\]
with Riemann symbol
\[
\left\{
\begin{array}{ccc}
0&1/4N&\infty\\
\hline
0&0&\alpha_1\\
0&1&\alpha_2\\
0&3/2&1/2\\
0&2&\alpha_3\\
0&3&\alpha_4\\
\end{array}
\right\}
\]
These operators have a {\em Yifan Yang pull-back} to $14$ special fourth order 
operators, called the {\em tilde}-operators $\widetilde{1},\widetilde{2},\ldots, 
\widetilde{14}$. These operators replace the more complicated
{\em hat-operators}  $\hat{i},i=1,2,\ldots,14$ that appeared in the list \cite{AESZ}
into which they can be transformed.\\
 
An explicit formula for the tilde-operators can be found in \cite{Alm3}:
\[
\theta^4-4Nx\left(2(\theta+\frac{1}{2})^4+\frac{1}{2}(\frac{7}{2}-\mu^2-\nu^2)(\theta+\frac{1}{2})^2+\frac{1}{16}-\frac{1}{4}(\mu^2+\frac{1}{4})(\nu^2+\frac{1}{4})\right)+\]
\[+(4N)^2x^2(\theta+1+\frac{\mu+\nu}{2})(\theta+1-\frac{\mu+\nu}{2})(\theta+1+\frac{\mu-\nu}{2})(\theta+1-\frac{\mu-\nu}{2})
\]
where $\mu:=\alpha_3-\frac{1}{2}=\frac{1}{2}-\alpha_2$, $\nu:=\alpha_4-\frac{1}{2}=\frac{1}{2}-\alpha_1$.
Its Riemann-symbol is of the form
\[
\left\{
\begin{array}{ccc}
0&1/4N&\infty\\
\hline
0&-1/2& 1-(\mu+\nu)/2 \\
0&0& 1+(\mu-\nu)/2 \\
0&1&1-(\mu-\nu)/2 \\
0&3/2&1+(\mu+\nu)/2 \\
\end{array}
\right\}
\]
A general formula for the homolorphic solution is also given in \cite{Alm3}.
The $2\times 2$ Wronskians for this operator, multiplied by
a factor $(1-4Nx)^{3/2}$, are solutions to the hypergeometric
fifth order equation.\\

The monodromy around the central point $1/4N$ is of order two, and the
matrix in the scaled Frobenius basis is of the form
\[
\left(
\begin{array}{cccc}
x & 0 & y& 0 \\
z & -x &  0 & y \\
t &  0 & -x  & 0 \\
0 & t & -z&  x \\
\end{array}
\right)
\]

where $x,y,t$ satisfy the relation $x^2+yt=1$.
The invariants $x,y,z,t$ can be expressed in terms of the reflection vector
(see section 4.3):
\[(a,0,c,d)=(h^3,0,c_2h/24,\chi \lambda), \;\;\; \lambda:=\frac{\zeta(3)}{(2\pi i)^3}\]
of the corresponding hypergeometric operator by the formulas
\[
x=\frac{\sqrt{a}}{24}+\frac{c}{\sqrt{a}},\;\;\;
y=\sqrt{a}\;\;\;z=\frac{2d-4a\lambda}{\sqrt{a}},\;\;\;
t=(1-x^2)/y.\]

\newpage
\centerline{\Large \bf \em The fourteen hypergeometric and the corresponding tilde operators}
\resizebox{\linewidth}{!}
{
$
\begin{array}{|c|cccc|c||c|c|}
\hline
\textup{Case}&&\textup{Exponents}&&&N&\textup{Monodromy data}&\textup{Number}\\
\hline
1& 1/5&2/5&3/4&4/5&5^5&5,0,25/12,-200\lambda&1.1\\
\widetilde{1}&4/5&9/10&11/10&6/5 & & 11\sqrt{5}/24, \sqrt{5}, -84\sqrt{5} \lambda -29\sqrt{5}/2880&2.39\\
\hline
2 & 1/10&3/10&7/10&9/10&2^85^5&1,0,17/12,-288 \lambda&1.2\\
\widetilde{2}&7/10&9/10&11/10&13/10& &35/24, 1,-580 \lambda, -649/576&2.40\\
\hline
3 &  1/2&1/2&1/2&1/2&2^8 &16,0,8/3,-128 \lambda&1.3\\
\widetilde{3} & 1& 1& 1& 1& &5/6 ,4,-80 \lambda, &\tilde 2.33\\
\hline
4 & 1/3&1/3&2/3&2/3&3^6 &9,0,9/4,-144 \lambda&1.4\\
\widetilde{4} &5/6&1&1&7/6& &7/8, 3, -108 \lambda, 5/64&2.41\\
\hline
5 &1/3&1/2&1/2&2/3 &2^43^3&12,0,5/2,-144 \lambda  &1.5\\
\widetilde{5} &11/12&11/12&13/12&13/12& &\sqrt{3}/2, 2\sqrt{3}, -56\sqrt{3}\lambda, \sqrt{3}/24&2.42\\
\hline
6 &1/4&1/2&1/2&3/4&2^{10}&8,0,7/3,-296 \lambda&1.6\\
\widetilde{6} &7/8&7/8&9/8&9/8& &2\sqrt{2}/3,2\sqrt{2},-96\sqrt{2}\lambda,\sqrt{2}/36&2.43\\
\hline
7 &1/8&3/8&5/8&7/8&2^{16}&2,0,11/6,-296 \lambda&1.7\\
\widetilde{7} &3/4&7/8&9/8&5/4& &23\sqrt{2}/24,\sqrt{2},-300\sqrt{2}\lambda,-241\sqrt{2}/576&2.44\\
\hline
8 &1/6&1/3&2/3&5/6&2^43^6&3,0,7/4,-204 \lambda&2.45\\
\widetilde{8} &3/4&11/12&13/12&5/4& &5\sqrt{3}/8,\sqrt{3},-140\sqrt{3}\lambda,-11\sqrt{3}/192 &2.45\\
\hline
 9 &1/12&5/12&7/12&11/12&2^{12}3^6&1,0,23/12, -484 \lambda&1.9\\
\widetilde{9} &3/4&5/6&7/6&5/4& &47/24,1,-972 \lambda,-1633/576&2.46\\
\hline
10 &1/4&1/4&3/4&3/4&2^{12}&4,0,5/3,-144 \lambda&1.10\\
\widetilde{10}&3/4&1&1&5/4 & &11/12,2,-152\lambda,23/288&2.47\\
\hline
11&1/4&1/3&2/3&3/4&2^63^3&6,0,2,-156 \lambda&1.11\\
\widetilde{11}&19/24&23/24&25/24&29/24 & &3\sqrt{6}/8,\sqrt{6},-56\sqrt{6}\lambda,5\sqrt{6}/192&2.48\\
\hline
12&1/5&1/4&3/4&5/6&2^{10}3^3&2,0,4/3,-156 \lambda&1.12\\
\widetilde{12}&17/24&23/24&25/24&31/24& &17\sqrt{2}/24,\sqrt{2},-160\sqrt{2}\lambda,-\sqrt{2}/576&2.49\\
\hline
 13&1/6&1/6&5/6&5/6&2^83^6&1,0,11/12,-120 \lambda&1.13\\
\widetilde{13}&2/3&1&1&4/3& &23/24,1,-244\lambda,47/576&2.50\\
\hline
14&1/6&1/2&1/2&5/6&2^83^3&4,0,13/6,-256 \lambda&1.14\\
\widetilde{14}&5/6&5/6&7/6&7/6& &7/6,2,-264\lambda,-13/72&2.51\\
\hline
\end{array}
$ 
}

\newpage
{\Large \bf \em List of Calabi-Yau operators of degree 2}\\

We give the parameters $a,b,c,d,e,f,\alpha,\beta,\gamma,\delta$ for the
Calabi-Yau operators written in the form (note the sign-change)
\[\theta^4-x(a\theta^4+b\theta^3+c\theta^2+d\theta+e)+f x^2(\theta+\alpha)(\theta+\beta)(\theta+\gamma)(\theta+\delta).\]

\resizebox{\columnwidth}{!}
{%
$
\begin{array}{|c||c|c|c|c|c|c||c||c|c|c|c|}
\hline
Number&a&b&c&d&e&f&\alpha&\beta&\gamma&\delta\\
\hline
2.1& 112& 224& 172& 60& 8   & - (2)^{11}& 1/2& 1/2& 3/2& 3/2\\
\hline
2.2 & 189& 378& 285& 96& 12  & - (2)^3(3)^6 & 1/3& 2/3& 4/3& 5/3\\
\hline
2.3 & 448& 896& 660& 212& 24  & - (2)^{15}  & 1/4& 3/4& 5/4& 7/4\\
\hline
2.4 & 3024& 6048& 4308& 1284& 120   & - (2)^{11}(3)^6 & 1/6& 5/6& 7/6& 11/6\\
\hline
2.5 & 176&  352& 268& 92& 12  & - (2)^8 &1/2& 1/2& 3/2& 3/2\\
\hline
2.6 & 4752& 9504& 6708& 1956& 180 & - (3)^6 & 1/3& 2/3& 4/3& 5/3\\
\hline
2.7 & 704& 1408& 1028& 324& 36  & - (2)^{12}  & 1/4& 3/4& 5/4& 7/4\\
\hline
2.8 &  4752& 9504& 6708& 1956& 180 & - (2)^8(3)^6 & 1/6& 5/6& 7/6& 11/6\\
\hline
2.9 &  160&320&248&88&12& (2)^8(3)^2& 1/2& 1/2& 3/2& 3/2\\
\hline
2.10&  270& 540& 411& 141& 18& (3)^8 & 1/3& 2/3& 4/3& 5/3\\
\hline
2.11&  640&1280&952&312&36&  (2)^{12}(3)^2 & 1/4& 3/4& 5/4& 7/4\\
\hline
2.12& 4320&8640&6216&1896&180&  (2)^8(3)^8& 1/6& 5/6& 7/6& 11/6\\
\hline
2.13& 192&384&304&112&16&  (2)^{13}  & 1/2& 1/2& 3/2& 3/2\\
\hline
2.14& 324&648&504&180&24&  (2)^5(3)^6&1/3&2/3& 4/3& 5/3\\
\hline
2.15& 768&1536&1168&400&48&  (2)^{17}  & 1/4& 3/4& 5/4& 7/4\\
\hline
2.16& 5184&10368&7632&2448&240&  (2)^{13}(3)^6 & 1/6& 5/6& 7/6& 11/6\\
\hline
2.17& 512&1024&832&320&48&  (2)^{16}  & 1/2& 1/2& 3/2& 3/2\\
\hline
2.18&  864&1728&1380&516&72&  (2)^8 (3)^6 & 1/3& 2/3& 4/3& 5/3\\
\hline
2.xx19&-144& -1152& -3200& -4096& -2048&  (2)^{20}  & 1/4& 3/4& 5/4& 7/4\\
\hline
2.19& 13824&27648&20928& 7104& 720&  (2)^{16}(3)^3 & 1/6& 5/6& 7/6& 11/6\\
\hline
2.20& 144&288&228&84&12&  (2)^8   (3)^3 & 1/2& 1/2& 3/2& 3/2\\
\hline
2.21& 243&486&378&135&18&  (3)^9 & 1/3& 2/3& 4/3& 5/3\\
\hline
2.22& 576& 1152& 876& 300&36&  (2)^{12}(3)^3 & 1/4& 3/4& 5/4& 7/4\\
\hline
2.23& 3888&7776&5724&1836&180&  (2)^8(3)^9 & 1/6& 5/6& 7/6& 11/6\\
\hline
2.24& 272&544&436&164&24&  (2)^{11} (3)^2 & 1/2& 1/2& 3/2& 3/2\\
\hline
2.25& 459&918&723&264&36&  (2)^3   (3)^8 & 1/3& 2/3& 4/3& 5/3\\
\hline
2.26& 1088& 2176&1676& 588& 72&  (2)^{15}(3)^2 & 1/4& 3/4& 5/4& 7/4\\
\hline
2.27& 7344&14688&10956&3612&360&  (2)^{11}    (3)^8 & 1/6& 5/6& 7/6& 11/6\\
\hline
2.28&864&1728&1416&552&84&  (2)^8   (3)^6 & 1/2& 1/2& 3/2& 3/2\\
\hline
2.29&1458& 2916& 2349& 891& 126&  (3)^{12}  & 1/3& 2/3& 4/3& 5/3\\
\hline
2.30& 23328& 46656& 35640& 12312& 1260&  (2)^8   (3)^{12}  & 1/6& 5/6& 7/6& 11/6\\
\hline
\end{array}
$
}

\resizebox{\columnwidth}{!}
{%
$
\begin{array}{|c|c|c|c|c|c|c|c|c|c|c|c|}
\hline
Number&a&b&c&d&e&f&\alpha&\beta&\gamma&\delta\\
\hline
2.31& 2048&4096&3392&1344&208&  (2)^{20}  & 1/2& 1/2& 3/2& 3/2\\
\hline
2.32& 3456&6912&5628&2172 &312&  (2)^{12} (3)^6 & 1/3& 2/3& 4/3& 5/3\\
\hline
2.33& 8192&16384&13056&4864&624&  (2)^{24}& 1/4& 3/4& 5/4& 7/4\\
\hline
2.34&  55296&110592&85440&30144&3120&  (2)^{20}(3)^6 & 1/6& 5/6& 7/6& 11/6\\
\hline
2.35& 13824&27648&23232&9408&1488&  (2)^{16} (3)^6 & 1/2& 1/2& 3/2& 3/2\\
\hline
2.36& 23328&46656&38556&15228 &2232&  (2)^8   (3)^{12}  & 1/3& 2/3& 4/3& 5/3\\
\hline
2.37&55296&110592&89472&34176&4464&  (2)^{20}    (3)^6 & 1/4&3/4& 5/4& 7/4\\
\hline
2.38& 373248&746496&585792&212544&22320&  (2)^{16}    (3)^{12}  & 1/6& 5/6& 7/6&11/6\\
\hline
2.39& 25000&50000&58750&33750&7380&  (2)^4 (5)^{10}  &4/5& 6/5& 9/10& 11/10\\
\hline
2.40&6400000&12800000&14880000&8480000&1824880&  (2)^{20}(5)^{10}& 7/10&9/10&11/10&13/10\\
\hline
2.41&5832&11664&13770&7938&1746&  (2)^4   (3)^{12}  & 1& 1& 5/6& 7/6\\
\hline
2.42& 3456&6912&8184&4728&1044&  (2)^{12} (3)^6 &11/12&11/12&13/12&13/12\\
\hline
2.43& 8192&16384&19328&11136&2448&  (2)^{24} & 7/8& 7/8& 9/8& 9/8\\
\hline
2.44& 524288&1048576&1224704&700416&151920&  (2)^{36}  & 3/4& 5/4& 7/8& 9/8\\
\hline
2.45&93312&186624&218376&125064&27180&  (2)^{12}(3)^{12}& 3/4& 5/4& 11/12&13/12\\
\hline
2.46& 23887872&47775744&55655424&31767552&6870384&  (2)^{28}(3)^{12}  & 3/4& 5/4& 5/6& 7/6\\
\hline
2.47&32768&65536&76800&44032&9584&  (2)^{28}  & 3/4& 1& 1& 5/4\\
\hline
2.48&  13824&27648&32520&18696&4092&  (2)^{16}(3)^6 &19/24&23/14&25/24&29/24\\
\hline
2.49&221184&442368&515712&294528&63600&  (2)^{24} (3)^{6} &17/24& 23/24& 25/24& 31/24\\
\hline
2.50&1492992&2985984&3462912&1969920&421488&  (2)^{20}(3)^{12}  & 1& 1& 2/3& 4/3\\
\hline
2.51& 55296&110592&129792&74496&16272&  (2)^{20}    (3)^{6} & 5/6& 5/6& 7/6& 7/6\\
\hline
2.52& 80&160&132&52&8&  (2)^{10}  & 1/2& 1& 1& 3/2\\
\hline
2.53& 136&272&210&74&10&  (2)^4 & 1/2& 1& 1& 3/2\\
\hline
2.54&56&112&94&38&6&  (2)^4   (3)^4 & 1/2& 1& 1& 3/2\\
\hline
2.55&96&192&152&56&8&  (2)^8 & 1/2& 1& 1& 3/2\\
\hline
2.56&72&144&114&42&6& - (2)^4   (3)^3 & 1/2& 1& 1& 3/2\\
\hline
2.57&88&176&150&62&10&  (2)^4   (5)^3 & 1/2& 1& 1& 3/2\\
\hline
2.xx58&128&256&224&96&16&(2)^{12} &1/2&1&1&3/2\\
\hline
2.yy58&2048&4096 &3200 &1152&144&(2)^{20}&1/4&3/4&5/4&7/4\\
\hline
2.58&216&432&390&174&30&  (2)^{4}   (3)^6 &1/2&1& 1& 3/2\\
\hline
2.59&3456&6912&6816&3360&624&  (2)^{12}    (3)^6 & 1/2& 1& 1& 3/2\\
\hline
2.60&48&96&76&28&4& - (2)^{10}  & 1/2& 3/4& 5/4& 3/2\\
\hline
2.61&104&208&162&58&8& - (2)^4   (3)^3 & 1/2& 2/3& 4/3& 3/2\\
\hline
2.62&65&130&105&40 &6&  (2)^6 & 3/4& 1& 1& 5/4\\
\hline
2.63&128&256&252&124&24&  (2)^{12}  & 3/4& 1& 1& 5/4\\
\hline
2.64&43&86&77&34&6&  (2)^4   (3) & 1& 1& 5/6& 7/6\\
\hline
2.65&112&224&188&76&12&  (2)^{10} (3)^6& 1/2& 3/4& 5/4& 3/2\\
\hline
2.66&864&1728&1416&552& 60&  (2)^4  (3)^4 & 1/6& 5/6& 7/6& 11/6\\
\hline
2.67&216& 396& 366& 168& 30&(2)^4(3^6)&2/3&2/3&7/6&7/6\\
\hline
2.68&512&896&788&340& 56& (2)^{16}& 1/2& 5/4& 5/8& 9/8\\
\hline
2.69&59&118&105&46&8&  (2)^5   (3)^3 & 1& 1& 2/3& 4/3\\
\hline
2.70&512&1408&1652&948&216&  (2)^{16}& 5/4& 3/2& 9/8& 13/8\\
\hline
\end{array}
$
}

\end{document}